\documentclass[12pt]{amsart}

\usepackage{epsf}
\usepackage{amsmath,amssymb}
\usepackage{colordvi}

\copyrightinfo{}{}

\usepackage[all]{xy}
\CompileMatrices
\newcommand{\xyinc}{\ar@{^{(}->}}

\newcommand{\map}[1]{\xrightarrow{#1}}

\numberwithin{equation}{section}
\theoremstyle{plain}
\newtheorem{thm}{Theorem}[section]
\newtheorem{prop}[thm]{Proposition}
\newtheorem{lemm}[thm]{Lemma}\newtheorem{coro}[thm]{Corollary}

\theoremstyle{definition}
\newtheorem{defi}[thm]{Definition}

\newtheorem{rem}[thm]{Remark} 

\newcommand{\sumsub}[1]{\sum_{\substack{#1}}} 

\newcommand{\ten}{\mbox{\hspace*{-.5pt}\raisebox{1pt}{${\scriptstyle \otimes}$}
\hspace*{-4pt}}}
\newcommand{\st}{\mathrm{st}}          
\newcommand{\id}{\mathit{id}}          

\newcommand{\Inv}{\mathrm{Inv}}

\newcommand{\Des}{\mathrm{Des}}
\newcommand{\GDes}{\mathrm{GDes}}             
\newcommand{\AGDes}{\overline{\mathrm{GDes}}} 

\newcommand{\setS}{\mathsf{S}}
\newcommand{\setT}{\mathsf{T}}
\newcommand{\setR}{\mathsf{R}}

\newcommand{\Q}{\mathbb{Q}}
\newcommand{\Z}{\mathbb{Z}}

\newcommand{\frakS}{\mathfrak{S}}

\newcommand{\calD}{\mathcal{D}}
\newcommand{\calF}{\mathcal{F}}

\newcommand{\calM}{\mathcal{M}}
\newcommand{\calQ}{\mathcal{Q}}
\newcommand{\calZ}{\mathcal{Z}}

\newcommand{\QSym}{\mathcal{Q}\mathit{Sym}}
\newcommand{\SSym}{\mathfrak{S}\mathit{Sym}}

\newcommand{\Sh}[1]{\mathfrak{S}^{#1}}

\DeclareMathOperator*{\Disjoint}{{\textstyle \coprod}}

\newcounter{FNC}[page]
\def\newfootnote#1{{\addtocounter{FNC}{2}$^\fnsymbol{FNC}$%
     \let\thefootnote\relax\footnotetext{$^\fnsymbol{FNC}$#1}}}

\title[Hopf algebra of permutations]{Structure of the Malvenuto-Reutenauer\\
       Hopf algebra of permutations}
\author{Marcelo Aguiar}
\address{Department of Mathematics\\
         Texas A\&M University\\
         College Station\\
         Texas \ 77843\\
         USA}
\email{maguiar@math.tamu.edu}
\urladdr{http://www.math.tamu.edu/$\sim$maguiar}

\author{Frank Sottile}
\address{Department of Mathematics\\
          Texas A\&M University\\
         College Station\\
         Texas \ 77843\\
         USA}
\email{sottile@math.tamu.edu}
\urladdr{http://www.math.tamu.edu/$\sim$sottile}

\thanks{Research of Sottile supported in part by NSF grant DMS-0070494}

\begin{document}

\begin{abstract}
 We analyze the structure of the Malvenuto-Reutenauer
 Hopf algebra $\SSym$ of permutations in detail.
 We give explicit formulas for its antipode, prove that it is a cofree
 coalgebra,  determine its  primitive elements and its coradical filtration,
 and show that it decomposes  as a crossed product over the Hopf algebra of
 quasi-symmetric functions. In addition, we describe the structure
 constants of the multiplication as a certain number of facets of the
 permutahedron.
 As a consequence we obtain a new interpretation of the product of
 monomial quasi-sym\-metric functions in terms of the facial structure of the
 cube.
 The Hopf algebra of Malvenuto and Reutenauer has a
 linear basis  indexed by permutations. Our results are obtained from a
 combinatorial description of the Hopf algebraic structure with respect to a new
 basis for this algebra, related to the original one via
 M\"obius inversion on the weak order on the symmetric groups.
 This is in analogy with the relationship between the monomial and fundamental
 bases of the algebra of quasi-sym\-metric functions.
 Our results reveal a close relationship between the structure of the
 Malvenuto-Reutenauer Hopf algebra and the weak order on the symmetric groups.
\end{abstract}

\maketitle
\contentsline {section}{\tocsection {}{}{Introduction}}{1}
\contentsline {section}{\tocsection {}{1}{Basic definitions and results}}{3}
\contentsline {section}{\tocsection {}{2}{The weak order on the symmetric group}}{7}
\contentsline {section}{\tocsection {}{3}{The coproduct of $\SSym$}}{14}
\contentsline {section}{\tocsection {}{4}{The product of $\SSym$}}{15}
\contentsline {section}{\tocsection {}{5}{The antipode of $\SSym$}}{18}
\contentsline {section}{\tocsection {}{6}{Cofreeness, primitive elements, and the coradical filtration of $\SSym$}}{23}
\contentsline {section}{\tocsection {}{7}{The descent map to quasi-symmetric functions}}{27}
\contentsline {section}{\tocsection {}{8}{$\SSym$ is a crossed product over $\QSym$}}{31}
\contentsline {section}{\tocsection {}{9}{Self-duality of $\SSym$ and applications}}{34}
\contentsline {section}{\tocsection {}{}{References}}{39}
\section*{Introduction}

Malvenuto~\cite{Malv} introduced the Hopf
algebra $\SSym$ of permutations, which has a linear
basis $\{\calF_u\mid u\in \frakS_n, n\geq0\}$ indexed by
permutations in all symmetric groups $\frakS_n$.
The Hopf algebra $\SSym$ is non-commutative, non-cocommutative,
self-dual, and graded.
Among its sub-, quotient-, and subquotient- Hopf algebras
are many algebras central to algebraic combinatorics.
These include the algebra of symmetric functions~\cite{Mac,St99},
Gessel's algebra $\QSym$ of quasi-symmetric functions~\cite{Ges}, the algebra
of non-commutative symmetric functions~\cite{GKal}, the Loday-Ronco
algebra of planar binary trees~\cite{LR98}, Stembridge's algebra 
of peaks~\cite{Stem97}, the Billera-Liu algebra of Eulerian
enumeration~\cite{BilLiu}, and others.
The structure of these combinatorial Hopf algebras with respect to certain
distinguished bases has been an important theme in algebraic combinatorics,
with applications to the combinatorial problems these algebras were created
to study.
Here, we obtain a detailed understanding of the
structure of $\SSym$, both in algebraic and combinatorial terms.

Our main tool is a new basis
$\{\calM_u\mid u\in \frakS_n, \ n\geq0\}$ for $\SSym$ related to the original
basis by M\"obius inversion on the weak order on the symmetric
groups. 
These bases $\{\calM_u\}$ and $\{\calF_u\}$ are analogous to the
monomial basis and the fundamental basis of $\QSym$,
which are related via M\"obius inversion on their index sets, the Boolean
posets $\calQ_n$. We refer to them as the monomial basis and the fundamental basis
of $\SSym$.

We give enumerative-combinatorial descriptions of the product, coproduct,
and antipode of $\SSym$ with respect to the monomial basis
$\{\calM_u\}$. 
In Section 3, we show that the coproduct is obtained by splitting a
permutation at certain special positions that we call global descents. 
Descents and global descents
are left adjoint and right adjoint to a natural map $\calQ_n\to\frakS_n$.
These results rely on some non-trivial
properties of the weak order developed in Section~\ref{S:bruhat}.

The product is studied in Section~\ref{S:product}.
The structure constants are non-negative integers with the following
geometric-combinatorial description.
The 1-skeleton of the permutahedron $\Pi_{n-1}$ is the Hasse diagram of the
weak order on $\frakS_n$.
The facets of the permutahedron are canonically isomorphic to products of lower
dimensional permutahedra. Say that
a facet isomorphic to $\Pi_{p-1}\times\Pi_{q-1}$ has type $(p,q)$.
Given $u\in\frakS_p$ and $v\in\frakS_q$, such a facet has a distinguished
vertex corresponding to $(u,v)$ under the canonical isomorphism.
Then, for $w\in\frakS_{p+q}$, the coefficient of
$\calM_w$ in $\calM_u\cdot\calM_v$ is the number of facets of the
permutahedron $\Pi_{p+q-1}$ of type $(p,q)$ with the property that the
distinguished vertex is below $w$ (in the weak order) and closer to $w$ than
any other vertex in the facet. 

In Section~\ref{S:antipode} we find explicit formulas for the antipode with
respect to both bases. 
The structure constants with respect to the monomial basis
have constant sign, as for $\QSym$.
The situation is more complicated for the fundamental basis, which may explain
why no such explicit formulas were previously known.
\smallskip

Elucidating the elementary structure of $\SSym$ with respect to the
monomial basis reveals further algebraic structures of
$\SSym$. In Section~\ref{S:cofree}, we show that $\SSym$ is a cofree
graded coalgebra. 
A consequence is that its coradical filtration (a filtration encapsulating 
the complexity of iterated coproducts) is the algebraic
counterpart of a filtration of the symmetric groups by certain lower order
ideals. In particular, we show that the space of primitive elements is spanned
by the set $\{\calM_u\mid u \text{ has no global descents}\}$. 
Cofreenes was shown by Poirier and Reutenauer~\cite{PR95} in dual form,
through the introduction of a different basis. 
The study of primitive elements was pursued from this point of view 
by Duchamp, Hivert, and Thibon~\cite{DHT01}. 
The generating function for
the graded space of primitive elements is 
 \[
   1-\frac{1}{\sum_{n\geq 0}n!\, x^n}\,.
 \]
Comtet essentially studied the combinatorics of global
descents~\cite[Exercise VI.14]{Co74}.
These results add an algebraic perspective to the pure combinatorics
he studied.

There is a well-known morphism of Hopf algebras $\SSym\to\QSym$ that maps
one fundamental basis onto the other, by associating to a permutation $u$ its
descent set $\Des(u)$. 
In Section~\ref{S:descentmap}, we describe this map on the monomial
bases and then derive a new geometric description for the product of
monomial quasi-symmetric functions in which the role of the permutahedron is
played by the cube.

In Section~\ref{S:crossed} we show that $\SSym$ decomposes as a crossed
product over $\QSym$. This construction from the theory of Hopf algebras
is a generalization of the notion of group extensions. We provide a
combinatorial description for the Hopf kernel of the map $\SSym\to\QSym$,
which is a subalgebra of $\SSym$.

We study the self-duality of $\SSym$ in Section~\ref{S:duality} and its
enumerative consequences.
For instance, a result of Foata and Sch\"utzenberger~\cite{FS78} on the
numbers  
 \[
    d(\setS,\setT)\ :=\ \#\{w\in\frakS_n\mid \Des(w) = \setS,\
                           \Des(w^{-1}) = \setT\}
 \]
follows directly from this self-duality and we obtain analogous
results for the numbers
 \[
   \theta(u,v)\ :=\ \#\{w\in\frakS_n\mid w\leq u,\ w^{-1}\leq v\}\,. 
\]
Most of the order-theoretic properties that underlie these algebraic
results are presented in Section~\ref{S:bruhat}. 
Central to these are the existence of two Galois connections (involving
descents and global descents) between the posets of permutations of $[n]$ and
of subsets of $[n{-}1]$, as well as the order properties of the decomposition
of $\frakS_n$ into cosets of $\frakS_p\times \frakS_q$.

We thank Swapneel Mahajan, who suggested a simplification to the
proof of Theorem~\ref{T:cop-monomial}, Nantel Bergeron, one of whose
questions motivated the results of Section~\ref{S:crossed}, and the referees
of an abridged version for helpful comments.

\section{Basic definitions and results}\label{S:basic}
We use only elementary properties of Hopf algebras, as given in 
the book~\cite{Mo93a}.
Our Hopf algebras $H$ will be graded connected Hopf algebras over $\Q$.
Thus the $\Q$-algebra $H$ is the direct sum 
$\bigoplus\{H_n\mid n=0,1,\ldots\}$ of its 
homogeneous components $H_n$, with $H_0=\Q$, the product and
coproduct respect the grading, and the counit is projection onto $H_0$.

Throughout, $n$ is a non-negative integer and $[n]$ denotes the set
$\{1,2,\ldots,n\}$.
A {\it composition} $\alpha$ of $n$ is a sequence
$\alpha=(\alpha_1,\ldots,\alpha_k)$ of positive integers with
$n=\alpha_1+\alpha_2+\cdots+\alpha_k$.
To a composition $\alpha$ of $n$, we associate the set $I(\alpha):=
\{\alpha_1,\alpha_1+\alpha_2,\ldots,\alpha_1+\cdots+\alpha_{k-1}\}$.
This gives a bijection between compositions of $n$ and subsets of $[n{-}1]$.
Compositions of $n$ are partially ordered by \emph{refinement}.
The cover relations are of the form
 \[
   (\alpha_1,\ldots,\alpha_i+\alpha_{i+1},\ldots,\alpha_k)\ \lessdot\ 
   (\alpha_1,\ldots,\alpha_k)\,.
 \]
Under the association $\alpha\leftrightarrow I(\alpha)$, refinement
corresponds to set inclusion, so we simply identify the poset of compositions
with the Boolean poset $\calQ_n$ of subsets of $[n{-}1]$.

Let $\frakS_n$ be the group of permutations of $[n]$.
We use one-line notation for permutations, writing $u=(u_1,u_2,\ldots,u_n)$
where $u_i=u(i)$.
Sometimes we may omit the commas and write $u=u_1\ldots u_n$.
A permutation $u$ has a \emph{descent} at a position $p$ if 
$u_p>u_{p+1}$.
An \emph{inversion} in a permutation $u\in\frakS_n$ is a pair of positions
$1\leq i<j\leq n$ with $u_i>u_j$. The set of descents and inversions are
denoted by $\Des(u)$ and $\Inv(u)$, respectively.
The length of a permutation $u$ is $\ell(u)=\#\Inv(u)$.

Given $p,q\geq 0$, we consider the product $\frakS_p\times\frakS_q$ to be a
subgroup of $\frakS_{p+q}$, where $\frakS_p$ permutes $[p]$ and $\frakS_q$
permutes $\{p{+}1,\ldots,p{+}q\}$.
For $u\in\frakS_p$ and $v\in\frakS_q$, write $u\times v$ for the permutation
in $\frakS_{p+q}$ corresponding to $(u,v)\in\frakS_p\times\frakS_q$
under this embedding.

More generally, given a subset $\setS=\{p_1<\cdots<p_k\}$ of $[n{-}1]$, 
we have the (standard) parabolic or Young subgroup
 \[
   \frakS_{\setS}\ :=\ \frakS_{p_1}\times\frakS_{p_2-p_1}
                        \times\cdots\times\frakS_{n-p_k}
          \ \subseteq\ \frakS_n\,.
 \]
The notation $\frakS_{\setS}$ suppresses the dependence on $n$, which will
either be understood or will be made explicit when this is used.

Lastly, we use $\Disjoint$ to denote disjoint union.

\subsection{The Hopf algebra of permutations of Malvenuto and Reutenauer}

Let $\SSym$ be the graded vector space over $\Q$ with basis 
$\Disjoint_{n\geq0}\frakS_n$, graded by $n$. 
This vector space has a graded Hopf algebra structure first considered in
Malvenuto's thesis~\cite[\S 5.2]{Malv} and in her work with
Reutenauer~\cite{MR95}. 
(In~\cite{DHT01}, it is called the algebra of free quasi-symmetric functions.)
Write $\calF_u$ for the basis element corresponding to $u\in\frakS_n$ 
for $n>0$ and $1$ for the basis element of degree $0$.

The product of two basis elements is obtained by shuffling the corresponding
permutations, as in the following example.
 \begin{eqnarray*}
  \calF_{\Blue{12}}\cdot\calF_{\Brown{312}} &=& \ \ \ 
    \calF_{\Blue{12}\Brown{534}}\,+\,
    \calF_{\Blue{1}\Brown{5}\Blue{2}\Brown{34}}\,
    +\,\calF_{\Blue{1}\Brown{53}\Blue{2}\Brown{4}}\,+\,
    \calF_{\Blue{1}\Brown{534}\Blue{2}}\,+\,
    \calF_{\Brown{5}\Blue{12}\Brown{34}}\\
    && +\,\calF_{\Brown{5}\Blue{1}\Brown{3}\Blue{2}\Brown{4}}\,+\,
    \calF_{\Brown{5}\Blue{1}\Brown{34}\Blue{2}}\,+\,
    \calF_{\Brown{53}\Blue{12}\Brown{4}}\,
    +\,\calF_{\Brown{53}\Blue{1}\Brown{4}\Blue{2}}\,+\,
    \calF_{\Brown{534}\Blue{12}}\,.
 \end{eqnarray*}

More precisely, for $p,q>0$, set
 \begin{align*}
   \Sh{(p,q)} \ :=&\
              \{\zeta\in \frakS_{p+q}\mid \zeta \mbox{ has at most one
              descent, at position $p$}\}\\
             = & \ \{\zeta\in \frakS_{p+q}\mid \zeta_1<\dotsb<\zeta_{p},\
                \zeta_{p+1}<\dotsb<\zeta_{n} \}\,.
 \end{align*}
This is the  collection of minimal (in length) representatives of left cosets
of $\frakS_p\times\frakS_q$ in $\frakS_{p+q}$. In the literature, they are
sometimes referred to as $(p,q)$-shuffles, but sometimes it is the inverses of
these permutations that carry that name. We will refer to them as
Grassmannian permutations. With these definitions, we describe the product.
For $u\in \frakS_p$ and $v\in \frakS_q$, set
 \begin{equation}\label{E:prod-fundamental}
  \calF_u \cdot \calF_v\ =\ \sum_{\zeta\in \Sh{(p,q)}}
                            \calF_{(u\times v)\cdot\zeta^{-1}}. 
\end{equation}
This endows $\SSym$ with the structure of a graded algebra with unit 1.

The algebra $\SSym$ is also a graded coalgebra with coproduct given by all
ways of splitting a permutation.
For a sequence $(a_1,\ldots,a_p)$ of distinct  integers, let
its {\it standard permutation}\newfootnote{Some authors call this 
flattening.} $\st(a_1,\ldots,a_p)\in\frakS_p$
be the permutation $u$ defined by  
\begin{equation}\label{E:st}
   u_i<u_j \iff a_i<a_j.
\end{equation}
For instance, $\st(625)=312$.
The coproduct $\Delta\colon\SSym\to\SSym\,\ten\SSym$ is defined by
\begin{equation}\label{E:cop-malvenuto}
  \Delta(\calF_u)\ =\ \sum_{p=0}^n \calF_{\st(u_1,\,\ldots,\,u_p)}\ten
                                   \calF_{\st(u_{p+1},\,\ldots,\,u_n)}\,,
\end{equation}
when $u\in\frakS_n$.
For instance, $\Delta(\calF_{42531})$ is 
\[
  1\ten\calF_{42531} +\calF_{1}\ten\calF_{2431} + \calF_{21}\ten\calF_{321} 
  +\calF_{213}\ten\calF_{21}+\calF_{3142}\ten\calF_{1} +\calF_{42531}\ten 1\,.
\]
$\SSym$ is a graded connected Hopf algebra~\cite[th\'eor\`eme 5.4]{Malv}.

We refer to the 
set $\{\calF_u\}$ as the {\em fundamental} basis of $\SSym$. 
The main goal of this paper is to obtain a detailed description of the Hopf
algebra structure of $\SSym$. To this end, the definition  of a
second basis for $\SSym$ (in \S~1.3) will prove crucial.

This Hopf algebra $\SSym$ of Malvenuto and Reutenauer has been an object of
recent interest~\cite{DHT, DHT01, Jol99, LR98, LR01,MR95, RP01, PR95, Re93}.
We remark that sometimes it is the dual Hopf algebra that is considered.
To compare results, one may use that $\SSym$ is
self-dual under the map $\calF_u\mapsto\calF_{u^{-1}}^*$, where
$\calF_{u^{-1}}^*$ is the element of the dual basis that is dual to 
$\calF_{u^{-1}}$.
We explore this further in Section~\ref{S:duality}.  

\subsection{The Hopf algebra of quasi-symmetric functions}\label{S:hopfquasi}
Basic references for quasi-symmetric functions are~\cite[9.4]{Re93}
and~\cite[Section 7.19]{St99}; however, everything we need will be reviewed
here.

The algebra $\QSym$ of quasi-symmetric functions is a subalgebra of the
algebra of formal power series in countably many variables
$x_1,x_2,\dotsc$.
It has a basis of {\em monomial quasi-symmetric functions} $M_\alpha$ indexed
by compositions $\alpha=(\alpha_1,\ldots,\alpha_k)$, where 
 \begin{equation*}
  M_\alpha\ :=\ \sum_{i_1<\dotsb<i_k}
      x_{i_1}^{\alpha_1}x_{i_2}^{\alpha_2}\cdots x_{i_k}^{\alpha_k}\,.
 \end{equation*}

The product of these monomial functions is given by quasi-shuffles of their
indices. 
A \emph{quasi-shuffle} of compositions $\alpha$ and $\beta$ is a shuffle of
the components of $\alpha$ and $\beta$, where in addition we may replace
any number of pairs of consecutive components $(\alpha_i,\beta_j)$ in the
shuffle by $\alpha_i+\beta_j$.
Then we have
 \begin{equation}\label{E:prodqsym}
   M_\alpha \cdot M_\beta\ =\ \sum_\gamma M_\gamma\,,
 \end{equation}
where the sum is over all quasi-shuffles $\gamma$ of the compositions $\alpha$
and $\beta$. For instance,
 \begin{equation}\label{E:prod-ex}
 M_{(2)} \cdot M_{(1,1)}  \ =\ 
   M_{(1,1,2)}+M_{(1,2,1)}+M_{(2,1,1)}+M_{(1,3)}+M_{(3,1)}\,.
 \end{equation}
The unit element is indexed by the empty composition $1=M_{(\ )}$.

Let $X$ and $Y$ be two countable ordered sets and $X\coprod Y$ its
disjoint union, totally ordered by $X<Y$.
Then $\Delta\colon f(X)\mapsto f(X\coprod Y)$ gives $\QSym$ the structure of a
coalgebra.
In terms of the monomial quasi-symmetric functions, we have
\begin{equation}\label{E:copqsym}
  \Delta\bigl(M_{(\alpha_1,\ldots,\alpha_k)}\bigr)\ =\,\ \sum_{p=0}^k
  M_{(\alpha_1,\ldots,\alpha_p)}\ten M_{(\alpha_{p+1},\ldots,\alpha_k)}\,.
\end{equation}
For instance, $\Delta(M_{(2,1)})=1\ten M_{(2,1)}+M_{(2)}\ten
M_{(1)}+M_{(2,1)}\ten 1$.

The algebra of quasi-symmetric functions was introduced by Gessel~\cite{Ges}.
Its Hopf algebra structure was
introduced by Malvenuto~\cite[Section 4.1]{Malv}. 
The description of the product in terms of quasi-shuffles can be
found in~\cite{Ho00} and is equivalent to~\cite[Lemma 3.3]{Eh96}.
A $q$-version of this construction appears in~\cite[Section 5]{TU96}
and~\cite{Ho00}.

The algebra $\QSym$ is a graded connected  Hopf algebra whose component in
degree $n$ is spanned by those $M_\alpha$ with $\alpha$ a composition of $n$.
Malvenuto~\cite[corollaire 4.20]{Malv} and 
Ehrenborg~\cite[Proposition 3.4]{Eh96} independently gave an explicit formula
for the antipode
 \begin{equation}\label{E:Q-antipode}
   S(M_\alpha)\ =\
   (-1)^{c(\alpha)}\sum_{\beta\leq\alpha}M_{\widetilde{\beta}}\,.
 \end{equation}
Here, $c(\alpha)$ is the number of components of $\alpha$, 
and if $\beta=(\beta_1,\beta_2,\ldots,\beta_t)$ then 
$\widetilde{\beta}$ is $\beta$ written in reverse order 
$(\beta_t,\ldots,\beta_2,\beta_1)$.

Gessel's {\em fundamental}  quasi-symmetric function $F_\alpha$ is defined by
 \[
    F_\alpha\ =\ \sum_{\alpha\leq\beta} M_\beta\,,
 \]
By M\"obius inversion, we have
 \[
    M_\alpha\ =\ \sum_{\alpha\leq\beta}(-1)^{c(\beta)-c(\alpha)}F_\beta\,.
 \]
Thus the set $\{F_\alpha\}$ forms another basis of $\QSym$. 

It is sometimes advantageous  to
index these monomial and fundamental quasi-symmetric functions by subsets of
$[n{-}1]$.
Accordingly, given a composition $\alpha$ of $n$ with $\setS=I(\alpha)$, we
define
 \[
    F_{\setS}\ :=\ F_\alpha \qquad\text{and}\qquad
    M_{\setS}\ :=\ M_\alpha\,.
 \]
The notation $F_{\setS}$ suppresses the dependence on $n$, which will be usually
understood from the context; otherwise it will be made explicit by writing
$F_{\setS,n}$.

In terms of power series,
 \begin{equation}   \label{E:QSym-def}
    F_{\setS}\ =\ \sumsub{i_1\leqslant i_2\leqslant\dotsb\leqslant i_n\\
       p\in\setS\Rightarrow i_p<i_{p+1}\rule{0pt}{8pt}}
       x_{i_1}x_{i_2}\cdots x_{i_n}   \,.
 \end{equation}
We mention that there is an analogous realization of the Malvenuto-Reutenauer
Hopf algebra as a subalgebra of an algebra of non-commutative power series, due
to Duchamp, Hivert, and Thibon.
To this end, one defines
 \begin{equation}   \label{E:SSym-def}
    \calF_u\ =\ \sumsub{i_1\leqslant i_2\leqslant\dotsb\leqslant i_n\\
      p\in\Des(u)\Rightarrow i_p<i_{p+1}\rule{0pt}{8pt}} 
      x_{i_{u^{-1}(1)}}x_{i_{u^{-1}(2)}}\cdots x_{i_{u^{-1}(n)}}   \,.
 \end{equation}
This is discussed in~\cite[Section 3.1]{DHT01}, in slightly different terms.
In this realization, the coproduct of $\SSym$ is induced by the ordinal sum of
commuting alphabets~\cite[Prop. 3.4]{DHT01}.

\subsection{The monomial basis of the Malvenuto-Reutenauer Hopf
            algebra}\label{S:monomialbasis} 

The descent set  of a permutation $u\in\frakS_n$ is the
subset of $[n{-}1]$ recording the descents of $u$
 \begin{equation}\label{E:defdescents}
   \Des(u)\ :=\ \{p\in[n{-}1]\mid u_p>u_{p+1}\}\,.
 \end{equation}
Thus $\Des(46\mspace{2mu}5\mspace{1mu}128\mspace{2mu}37)=\{2,3,6\}$.
Malvenuto~\cite[th\'eor\`emes 5.12, 5.13, and 5.18]{Malv} shows that 
there is a morphism of Hopf algebras
 \begin{equation} \label{E:descentmap}
   \begin{array}{rcrcl}
     \calD &:& \SSym&\longrightarrow& \QSym\\
           & &\calF_u&\longmapsto& F_{\Des(u)}\rule{0pt}{14pt}
   \end{array}
 \end{equation}
(This is equivalent to Theorem 3.3 in~\cite{MR95}.)
This explains our name and notation for the fundamental basis of $\SSym$.
This map extends to power series, where it is simply the abelianization: there is a
commutative diagram
\[\xymatrix{
 {\quad \SSym\quad }\xyinc[r]\ar@{>>}_{\calD}[d] & 
  \quad k\langle x_1,x_2,\ldots\rangle\quad \ar@{>>}^{\mathit{ab}}[d] \\ 
  {\quad \QSym\quad }\xyinc[r] &  \quad k[x_1,x_2,\ldots]\quad  }
\]
This is evident from~\eqref{E:QSym-def} and~\eqref{E:SSym-def}.
It is easy to see, however, that $\calD$ is {\it not} the 
abelianization of $\SSym$.

In analogy to the basis of monomial quasi-symmetric functions,
we define a new \emph{monomial} basis $\{\calM_u\}$ for the
Malvenuto-Reutenauer Hopf algebra. 
For each $n\geq 0$ and $u\in \frakS_n$, let
 \begin{equation}\label{E:def-monomial}
   \calM_u\ :=\ \sum_{u\leq v} \mu_{\frakS_n}(u,v)\cdot \calF_v\,,
 \end{equation}
where $u\leq v$ in the weak order in $\frakS_n$ 
(described in Section~\ref{S:bruhat}) and $\mu_{\frakS_n}$ 
is the M\"obius function of this partial order. 
By M\"obius inversion, 
 \begin{equation}\label{E:fun-mon}
   \calF_u\ :=\ \sum_{u\leq v} \calM_v\,,
 \end{equation}
so these elements $\calM_u$ indeed form a basis of $\SSym$. 
For instance,
 \[
   \calM_{4123}\ =\ \calF_{4123}-\calF_{4132}-\calF_{4213}+\calF_{4321}\,.
 \]
We will show that $\calM_u$ maps either to $M_{\Des(u)}$ or to 0 under the map
$\calD\colon\SSym\to\QSym$.

\section{The weak order on the symmetric group} \label{S:bruhat}

Let $\Inv(u)$ be the set of inversions of a permutation $u\in\frakS_n$,
 \[
   \Inv(u)\ :=\ \{(i,j)\in [n]\times[n]\mid i<j \text{ and }u_i>u_j\}\,.
 \]
The inversion set determines the permutation.
Given $u$ and $v\in\frakS_n$, we write $u\leq v$ if 
$\Inv(u)\subseteq\Inv(v)$. 
This defines the {\em left weak order} on $\frakS_n$.
\begin{figure}[htb]
  $$  \epsfxsize=3in\epsfbox{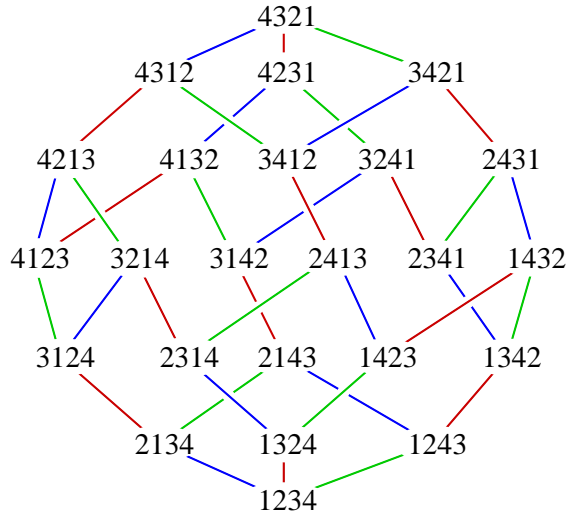}$$
  \caption{The weak order on $\frakS_4$\label{F:S4}}
\end{figure}
Figure~\ref{F:S4} shows the (left) weak order on $\frakS_4$.
The weak order has another characterization
 \[
   u\leq v\iff \exists w\in\frakS_n 
      \text{ such that }v=wu \text{ and } \ell(v)=\ell(w)+\ell(u)\,,
 \]
where $\ell(u)$ is the number of inversions of $u$.
The cover relations $u\lessdot v$ occur when $w$ is an
adjacent transposition. 
Thus, $u\lessdot v$ precisely when $v$ is obtained from $u$ by transposing a
pair of consecutive values of $u$; a pair $(u_i,u_j)$ such that $i<j$ and
$u_j=u_i+1$. 
The identity permutation $1_n$ is the minimum element of $\frakS_n$ and  
$\omega_n=(n,\ldots,2,1)$ is the maximum. 

This weak order is a lattice~\cite{GR60a}, whose structure we describe.
First, a set $J$ is the inversion set  of a permutation in $\frakS_n$ if and
only if both $J$ and its complement $\Inv(\omega_n) - J$ are 
{\em transitively closed} ($(i,j)\in J$ and $(j,k)\in J$ imply
$(i,k)\in J$, and the same for its complement).
The join (least upper bound) of two permutations $u$ and $v\in\frakS_n$ is the
permutation $u\vee v$ whose inversion set is 
the transitive closure of the union of the inversion sets of $u$ and $v$
 \begin{equation}\label{E:defjoin}
    \{(i,j)\mid\exists\text{ chain $i=k_0<\cdots<k_s=j$ s.t. 
              $\forall r,\ (k_{r-1},k_r)\in\Inv(u)\cup\Inv(v)$}\}.
 \end{equation}
Similarly, the meet (greatest lower bound) of $u$ and $v$ is the permutation
$u\wedge v$ whose inversion set is
 \begin{equation}\label{E:defmeet}
   \{(i,j)\mid\forall\text{ chains $i=k_0<\cdots<k_s=j,\ \exists r$ s.t.
    $(k_{r-1},k_r)\in\Inv(u)\cap\Inv(v)$}\}.
 \end{equation}

The M\"obius function of the weak order takes values in $\{-1,0,1\}$. 
Explicit descriptions  can be found in~\cite[Corollary 3]{Bj} 
or~\cite[Theorem 1.2]{Ede}.  
We will not need that description, but will use several basic facts on
the weak order that we develop here.

\subsection{Grassmannian permutations and the weak order}\label{S:shuffleweak}
In Section~\ref{S:basic}, we defined $\Sh{(p,q)}$ to be the set of minimal
(in length) representatives of (left) cosets of
$\frakS_p\times\frakS_q$ in the symmetric group $\frakS_{p+q}$.
Thus the map
 \[
   \begin{array}{rcccc}
    \lambda&:&\Sh{(p,q)}\times\frakS_p\times\frakS_q
              &\longrightarrow&\frakS_{p+q}\\ 
           & &(\zeta,u,v)&\longmapsto& \zeta\cdot(u\times v)\rule{0pt}{14pt}
    \end{array}
 \]
is a bijection.
We leave the following description of the inverse to the reader.

\begin{lemm}\label{L:stdecom}
 Let $w\in\frakS_{p+q}$, and set 
 $ \zeta :=\ w \cdot \bigl( \st(w_1,\ldots,w_p) \times 
                           \st(w_{p+1},\ldots,w_{p+q})\bigr)^{-1}$.
 Then $\zeta \in \Sh{(p,q)}$ and 
 $\lambda^{-1}(w)=(\zeta,\, \st(w_1,\ldots,w_p),\,
                          \st(w_{p+1},\ldots,w_{p+q}))$.
\end{lemm}

We describe the order-theoretic properties of this decomposition into cosets.The first step is to characterize the inversion sets of Grassmannian
permutations.
A subset $J$ of $[p]\times [q]$ is {\em cornered} if
$(h,k)\in J$ implies that $(i,j)\in J$ whenever 
$1\leq i\leq h$ and $1\leq j\leq k$.
The reason for this definition is that a 
set $I$ is the inversion set of a
Grassmannian permutation $\zeta\in \Sh{(p,q)}$ if and only if 
 \begin{equation}\label{L:invshuffle} 
   \begin{array}{rl}
    (i)& I\subseteq\{1,\ldots,p\}\times\{p{+}1,\ldots,p{+}q\}\,,
                 \text{ and}\\
    (ii)&\text{the shifted set }\{(p{+}1{-}i,j{-}p)\mid(i,j)\in I\}
                 \subseteq [p]\times[q] \text{ is cornered}\,.
                 \rule{0pt}{14pt}
   \end{array}
 \end{equation}
Given an arbitrary subset $J$ of $[p]\times[q]$, let $cr(J)$ denote the smallest
cornered subset containing $J$.
Denote the obvious action of $(u,v)\in\frakS_p\times \frakS_q$ on a subset $J$
of $[p]\times[q]$ by $(u,v)(J)$.

\begin{lemm}\label{L:corners} 
 Let $J$ be a cornered subset of $[p]\times[q]$ and $u\in\frakS_p$ and
 $v\in\frakS_q$ any permutations. 
 Then 
 \[
    J\ \subseteq\ \mathrm{cr}\bigl((u,v)(J)\bigr)\,.
 \]
\end{lemm}

\begin{proof}
 Let $(i,j)\in J$. 
 The set $\{u(h)\mid 1\leq h\leq i\}$ has $i$ elements. 
 Hence there is a number $h$ such that $1\leq h\leq i$ and $u(h)\geq i$. 
 Similarly there is number $k$ such that $1\leq k\leq j$ and $v(k)\geq j$.
 Since $J$ is cornered, $(h,k)\in J$. 
 Hence $(u(h),v(k))\in (u,v)(J)$. 
 By construction, $i\leq u(h)$ and $j\leq v(k)$, 
 so $(i,j)\in \mathrm{cr}\left((u,v)(J)\right)$, as needed.
\end{proof}

Denote the diagonal action of $w\in\frakS_n$  on a subset $I$
of $[n]\times [n]$ by $w(I)$.  Suppose $w=u\times v\in\frakS_p\times \frakS_q$
and $I\subseteq\{1,\ldots,p\}\times\{p{+}1,\ldots,p{+}q\}$. Let $J$ be the
result of shifting $I$, as  in~\eqref{L:invshuffle}({\it ii}). It is easy to see
that the result of shifting $(u\times v)(I)$ is $(\tilde{u},v)(J)$, where
$\tilde{u}(i)=p{+}1{-}u(p{+}1{-}i)$.

\begin{coro}\label{C:invshuffle}
 Let $\zeta$ and $\zeta'\in \Sh{(p,q)}$ be Grassmannian permutations, 
 and $u\in\frakS_p$ and $v\in\frakS_q$ be permutations. 
 If\/ $(u\times v)\left(\Inv(\zeta)\right)\subseteq\Inv(\zeta')$ then
 $\zeta\leq\zeta'$.
\end{coro}

\begin{proof} 
 We show that $\Inv(\zeta)\subseteq\Inv(\zeta')$. 
 Let $J$ and  $J'$ be the corresponding shifted sets. 
 According to the previous discussion and the hypothesis, 
 $(\tilde{u},v)(J)\subseteq J'$. Hence
 cr$\left((\tilde{u},v)(J)\right)\subseteq J'$,  since $J'$ is cornered.
 By Lemma~\ref{L:corners}, $J\subseteq
 \mathrm{cr}\left((\tilde{u},v)(J)\right)$,  so $J\subseteq J'$.
 This implies the inclusion of inversion sets, as needed. 
\end{proof}

The following lemma is straightforward.

\begin{lemm}\label{L:invlambda}
 Let $\zeta\in \Sh{(p,q)}$, $u\in\frakS_p$, $v\in\frakS_q$ and 
 $w:=\zeta\cdot(u\times v)\in\frakS_{p+q}$. 
 There is a decomposition of\/ $\Inv(w)$ into disjoint subsets 
 \[
   \Inv(w)\ =\ \Inv(u)\ \Disjoint\ \Bigl((p,p)+\Inv(v)\Bigr)
                  \ \Disjoint\ 
               (u^{-1}\times v^{-1})\Bigl(\Inv(\zeta)\Bigr)\,.
 \]
\end{lemm}

We deduce some order-theoretic properties of the decomposition into left
cosets. 
Define $\zeta_{p,q}$ to be the permutation of maximal length in
$\Sh{(p,q)}$, so that
 \[
   \zeta_{p,q}\ :=\ (q{+}1,\,q{+}2,\,\ldots,\,q{+}p,\ 1,\,2,\,\ldots,\,q)\,.
 \]

\begin{prop}\label{P:shuffles}
 Let $\lambda\colon\Sh{(p,q)}\times\frakS_p\times\frakS_q\to\frakS_{p+q}$
 be the bijection
 \[
   \lambda(\zeta,u,v)\ =\ \zeta\cdot(u\times v)\,.
 \]
 Then
 \begin{itemize}
  \item[({\it i})] $\lambda^{-1}$ is order preserving.
    That is,
   \[
     \zeta\cdot(u\times v)\ \leq\ \zeta'\cdot(u'\times v')\quad  
       \Longrightarrow\quad \zeta\ \leq\ \zeta',\ \ u\ \leq\  u',\
             \text{ and }\ v\ \leq\ v'\,.
    \]
  \item[({\it ii})] $\lambda$ is order preserving when restricted to any of the
             following sets
 \[
   \{\zeta_{p,q}\}\times\frakS_p\times\frakS_q,\ \
   \{1_{p+q}\}\times\frakS_p\times\frakS_q,\ \text{ or }\
    \Sh{(p,q)}\times\{(u,v)\}\,,
 \]
  for any $u\in\frakS_p$, $v\in\frakS_q$.
\end{itemize}
\end{prop}

\begin{proof}
 Let $w=\zeta\cdot(u\times v)$ and $w'=\zeta'\cdot(u'\times v')$.
 Suppose $w\leq w'$, so that $\Inv(w)\subseteq\Inv(w')$.
 By Lemma~\ref{L:invlambda}, we have $\Inv(u)\subseteq\Inv(u')$,
 $\Inv(v)\subseteq\Inv(v')$, and
 $(u''\times v'')\left(\Inv(\zeta)\right)\subseteq\Inv(\zeta')$,
 where $u'':=u'u^{-1}$ and $v'':=v'v^{-1}$.
 Therefore, $u\leq u'$, $v\leq v'$, and
 by Corollary~\ref{C:invshuffle}, $\zeta\leq\zeta'$.
 This proves ({\it i}).

 Statement ({\it ii}) follows by a similar application of
 Lemma~\ref{L:invlambda},
 noting that $\Inv(\zeta_{p,q})=\{1,\ldots,p\}\times\{p+1,\ldots,n\}$ and
 $\Inv(1_{p+q})=\emptyset$ are invariant under any permutation
 in $\frakS_p\times\frakS_q$.
\end{proof}

Since Grassmannian permutations in $\frakS^{(p,q)}$ are left coset
representatives of $\frakS_p\times\frakS_q$ in $\frakS_{p+q}$, their inverses
are right coset representatives.
We discuss order-theoretic properties of this decomposition into right
cosets.

Given a subset $J$ of $[n]\times[n]$, let
 \[
   \widetilde{J}\ =\ \{(j,i)\mid(i,j)\in J\}\,.
 \]
We have the following key observation about the diagonal action of $\frakS_n$
on subsets of $[n]\times[n]$.
\begin{lemm}\label{L:inversioninv}
 For any $u\in\frakS_n$, we have
 $u\bigl(\widetilde{\Inv(u)}\bigr) =\Inv(u^{-1})$.
\end{lemm}

\noindent{\it Proof.}
 Note that $u^{-1}(u_i)=i$.
 Thus $\Inv(u^{-1})=\{u_h<u_k\mid h>k\}$.
 Then
 \[
   u^{-1}\bigl(\Inv(u^{-1})\bigr)\ =\
   \{ (h,k)\mid k<h \text{ and }u_k>u_h\}\ =\
   \widetilde{\Inv(u)}\,. \eqno{\Box}\raisebox{-15pt}{\rule{0pt}{0pt}}
 \]

\begin{prop}\label{P:rightshuffling}
 Fix $\zeta\in\Sh{(p,q)}$ and consider the map
 $\rho_\zeta:\frakS_p\times\frakS_q\to\frakS_{p+q}$ given by
 \[
   \rho_\zeta(u,v)\ =\ (u\times v)\cdot\zeta^{-1}\,.
 \]
 Then $\rho_\zeta$ is a convex embedding in the sense that
 \begin{itemize}
   \item[(a)] $\rho_\zeta$ is injective;
   \item[(b)] $\rho_\zeta$ is order-preserving:
              $u\leq u' \text{ and }v\leq v' \iff
              \rho_\zeta(u, v)\leq \rho_\zeta(u',v')$;
   \item[(c)] $\rho_\zeta$ is convex:
              If $\rho_\zeta(u,v)\leq w\leq\rho_\zeta(u',v')$, for some
              $u,u'\in\frakS_p$ and $v,v'\in\frakS_q$, then there are
              $u''\in\frakS_p$ and $v''\in\frakS_q$ with
              $w=\rho_\zeta(u'',v'')$.
 \end{itemize}
 It follows that
 \begin{itemize}
   \item[(d)] $\rho_\zeta$ preserves meets and joins.
 \end{itemize}
\end{prop}
\begin{proof}
 Assertion (a) is immediate.
 Set $w:=(u\times v)\cdot\zeta^{-1}=\rho_\zeta(u,v)$.
 Then $w^{-1}=\zeta\cdot(u^{-1}\times v^{-1})$.
 By Lemmas~\ref{L:invlambda} and~\ref{L:inversioninv}, we have
 \begin{eqnarray*}
   \Inv(w)&=& w^{-1}\bigl(\widetilde{\Inv(w^{-1})}\bigr)\\
          &=&\zeta\cdot(u^{-1}\times v^{-1})\left(
              \widetilde{\Inv(u^{-1})}\cup
              \left((p,p)+\widetilde{\Inv(v^{-1})}\right)\cup
              (u\times v)\left(\widetilde{\Inv(\zeta)}\right)\right)\\
          &=& \zeta\left(\Inv(u)\cup \bigl((p,p)+\Inv(v)\bigr)\cup
                 \widetilde{\Inv(\zeta)}\right)\ .
 \end{eqnarray*}
 Assertion (b) follows from this and the characterization of the weak order in
 terms of inversion sets.

 For (c), decompose $w=(u''\times v'')\cdot{\xi}^{-1}$.
 By assumption,
 \[
   \zeta\Bigl(\widetilde{\Inv(\zeta)}\Bigr)\ \subseteq\
   \xi\Bigl(\widetilde{\Inv(\xi)}\Bigr)\ \subseteq\
   \zeta\Bigl(\widetilde{\Inv(\zeta)}\Bigr)\,.
 \]
 Then $\zeta=\xi$ by Lemma \ref{L:inversioninv}, so
 $w=\rho_\zeta(u'',v'')$ as needed.
\end{proof}

\subsection{Cosets of parabolic subgroups and the weak order}\label{S:multiweak}
Write a subset $\setS$ of $[n{-}1]$ as $\setS=\{p_1<\cdots<p_k\}$.
In Section~\ref{S:basic}, we defined the parabolic subgroup
 \[
   \frakS_\setS\ =\ \frakS_{p_1}\times\frakS_{p_2-p_1}\times\cdots\times
                     \frakS_{n-p_k}\ \subseteq\ \frakS_n\,.
 \]
Let $\Sh{\setS}$ be the set of minimal (in length) representatives of
left cosets of $\frakS_{\setS}$ in $\frakS_n$, 
\[
  \Sh{\setS}\ =\ \{\zeta\in\frakS_n \mid \Des(\zeta)\subseteq\setS\}\,.
\]
Grassmannian permutations are
the special case $\Sh{(p,n-p)}=\Sh{\{p\}}$.

Let $\zeta_{\setS}$ be the permutation of maximal length in $\Sh{\setS}$,
 \begin{equation}\label{E:maxshuffle}
   \zeta_\setS\ :=\
     (n{-}p_1{+}1,\ldots,n,\ n{-}p_2{+}1,\ldots,n{-}p_1,\ \dotsc,
        \ 1,\ldots,n{-}p_k)\,.
 \end{equation}
We record the following facts about these coset
representatives.

\begin{lemm}\label{L:galois}
 $\Sh{\setS}$ is an interval in the weak order of\/ $\frakS_n$.
 The minimum element is the identity $1_n$ and the maximum is $\zeta_\setS$.
\end{lemm}

Our proofs rely upon the following  basic fact.
Suppose $p,q$ are positive integers and $\setT$ is a subset of $[p{-}1]$.
Define the subset $\setS$ of $[p{+}q{-}1]$ to be $\setT\cup\{p\}$.
Then $(\zeta,\zeta')\mapsto \zeta\cdot(\zeta'\times 1_q)$ defines a bijection
 \begin{equation}\label{E:highershuffles}
   \Sh{(p,q)}\times \Sh{\setT}\ \longrightarrow\ \Sh{\setS}
 \end{equation}
The maximum elements are preserved under this map
 \begin{equation}\label{E:maxshuffles}
  \zeta_{p,q}\cdot(\zeta_\setT\times 1_q)\ =\ \zeta_\setS\,.
 \end{equation}

The analog of Proposition~\ref{P:shuffles} for this decomposition of
$\frakS_n$ into left cosets of $\frakS_\setS$ follows from
Proposition~\ref{P:shuffles}  by induction using~\eqref{E:highershuffles}
and~\eqref{E:maxshuffles}.

\begin{prop}\label{P:multishuffles}
 Suppose $\setS$ is a subset of $[n-1]$.
 Let $\lambda\colon\Sh{\setS}\times \frakS_\setS\to\frakS_n$ be the
 bijection
 \[
   \lambda(\zeta,u)\ =\ \zeta\cdot u\,.
 \]
Then $\lambda^{-1}$ is order preserving, while $\lambda$ is order
preserving when restricted to any of the
following sets
 \[
   \{\zeta_\setS\}\times\frakS_\setS,\ \
   \{1_{n}\}\times\frakS_\setS,\ \text{ or }\
   \Sh{\setS}\times\{u\}, \text{ for any $u\in\frakS_\setS$}\,.
 \]
\end{prop}

We state the analog of Proposition~\ref{P:rightshuffling}.

\begin{prop}\label{P:multirightshuffling}
 Let $\setS$ be a subset of\/ $[n{-}1]$.
 Fix $\zeta\in \Sh{\setS}$ and consider the map
 $\rho_\zeta\colon\frakS_\setS\to\frakS_n$ given by
 \[
    \rho_\zeta(u)\ =\ u\cdot\zeta^{-1}\,.
 \]
 Then $\rho_\zeta$ is a convex embedding.
 In particular, it preserves meets
 and joins.
\end{prop}

\subsection{Descents} \label{S:descents}

Let $\calQ_n$ denote the Boolean poset of subsets of $[n{-}1]$,
which we identify with the poset of compositions of $n$.
We have the descent map $\Des:\frakS_n\to\calQ_n$ given by
$u\mapsto \Des(u)$, the descent set~\eqref{E:defdescents} of $u$.
Let $Z:\calQ_n\to\frakS_n$ be the map defined by
$\setS\mapsto\zeta_\setS$,
the maximum left coset representative of $\frakS_\setS$ as
in~\eqref{E:maxshuffle}.

A \emph{Galois connection} between posets $P$ and $Q$ is a pair $(f,g)$ of
order preserving maps $f\colon P\to Q$ and $g\colon Q\to P$ such that for any
$x\in P$ and $y\in Q$,
 \begin{equation}\label{E:galois}
   f(x)\ \leq\ y \ \iff\  x\ \leq\ g(y)\,.
 \end{equation}

\begin{prop}\label{P:galois}
 The pair of maps $(\Des,Z):\frakS_n\rightleftarrows\calQ_n$
 is a Galois connection.
\end{prop}

\begin{proof}
 We verify that
 \begin{itemize}
   \item[(a)] $\Des\colon \frakS_n\to\calQ_n$ is order preserving;
   \item[(b)] $Z\colon \calQ_n\to\frakS_n$ is order preserving;
   \item[(c)] $\Des\circ Z=id_{\calQ_n}$;
   \item[(d)] $Z(\setS)=\max\{u\in\frakS_n\mid \Des(u)=\setS\}$.
 \end{itemize}

 First of all, the map $\Des$ is order preserving simply because $p$ is a
 descent  of $u$ if and only if $(p,p{+}1)\in\Inv(u)$.
 This is (a).
 The remaining assertions follow immediately from
 \begin{align*}
   \zeta_\setS\ &=\ \max\{u\in\frakS_n\mid \Des(u)\subseteq \setS\}\\
                 &=\ \max\{u\in\frakS_n\mid \Des(u) = \setS\}\,,
 \end{align*}
 which we know from Lemma~\ref{L:galois}.

 Condition~\eqref{E:galois} follows formally.
 In fact, suppose $\setT=\Des(u)\subseteq\setS$.
 Then by (d), $u\leq Z(\setT)$, and by (b),
 $Z(\setT)\leq Z(\setS)$, so $u\leq Z(\setS)$.
 Conversely, suppose $u\leq Z(\setS)$.
 Then by (a) and (c), $\Des(u)\subseteq \Des(Z(\setS))=\setS$.
\end{proof}

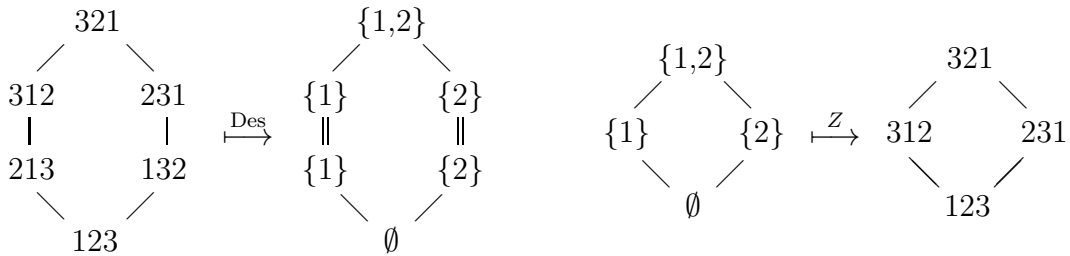
\begin{figure}[htb]
\[
  \begin{picture}(180,92)
   \put(0,0){
    \begin{picture}(70,92)
                        \put(25,83){321}
    \put(11,68){\line(1,1){10}}  \put(55,68){\line(-1,1){10}}
    \put( 0,55){312}             \put(50,55){231}
    \put( 8,40){\line(0,1){10}}  \put(60,40){\line(0,1){10}}
    \put( 0,27){213}             \put(50,27){132}
    \put(21,12){\line(-1,1){10}} \put(45,12){\line(1,1){10}}
                        \put(24,0){123}
    \end{picture}
   }

   \put(85,40){$\stackrel{\Des}{\longmapsto}$}

   \put(113,0){
    \begin{picture}(70,92)(1,0)
                        \put(19.5,83){\{1,2\}}
    \put(11,68){\line(1,1){10}}  \put(55,68){\line(-1,1){10}}
    \put(-1,55){\{1\}}             \put(50,55){\{2\}}
    \put( 7,40){\line(0,1){10}}  \put(58,40){\line(0,1){10}}
    \put( 9,40){\line(0,1){10}}  \put(60,40){\line(0,1){10}}
    \put(-1,27){\{1\}}             \put(50,27){\{2\}}
    \put(24,11){\line(-1,1){10}} \put(42,11){\line(1,1){10}}
                        \put(30,0){$\emptyset$}
    \end{picture}
   }
   \end{picture}
   \qquad\qquad
  \begin{picture}(170,92)
   \put(0,0){
    \begin{picture}(70,92)(1,0)
                        \put(19.5,69){\{1,2\}}
    \put(11,54){\line(1,1){10}}  \put(55,54){\line(-1,1){10}}
    \put(-1,41){\{1\}}             \put(50,41){\{2\}}
    \put(24,24){\line(-1,1){10}} \put(42,24){\line(1,1){10}}
                        \put(30,13){$\emptyset$}
    \end{picture}
   }

   \put(80,40){$\stackrel{Z}{\longmapsto}$}

   \put(103,0){
    \begin{picture}(70,92)
                        \put(25,69){321}
    \put(11,54){\line(1,1){10}}  \put(55,54){\line(-1,1){10}}
    \put( 2,41){312}             \put(53,41){231}
    \put(22,24){\line(-1,1){11}} \put(43,24){\line(1,1){11}}
                        \put(24,13){123}
    \end{picture}
   }
  \end{picture}
\]
 \caption{The Galois connection $\frakS_3\rightleftarrows\calQ_3$\label{F:galois}}
\end{figure}

This Galois connection is why the monomial basis  of
$\SSym$ is truly analogous to that of $\QSym$, and explains why we consider
the weak order rather than any other order on $\frakS_n$.
The connection between the monomial bases of these two algebras
will be elucidated in Theorem~\ref{T:map-monomial} using the previous result.

\subsection{Global descents}\label{S:global}

\begin{defi}\label{D:global}
 A permutation $u\in\frakS_n$ has a \emph{global descent} at a position
 $p\in[n{-}1]$ if
  \[ \forall\ i\leq p \text{ and }j\geq p{+}1\,,\ u_i>u_j\,.\]
Equivalently, if $\{u_1,\ldots,u_p\}=\{n,n{-}1,\ldots,n{-}p{+}1\}$.
 Let $\GDes(u)\subseteq[n{-}1]$ be the set of global descents of $u$.
 Note that $\GDes(u)\subseteq \Des(u)$, but these are not equal in
 general.
\end{defi}

In Section~\ref{S:descents} we showed that the descent map
$\Des\colon\frakS_n\to\calQ_n$ is left
adjoint to the map $Z\colon\calQ_n\to\frakS_n$, in the sense that the pair
$(\Des,Z)$ forms a Galois connection, as in Proposition~\ref{P:galois}.
That is, for $u\in\frakS_n$ and $\setS\in\calQ_n$,
 \begin{equation}\label{E:galoisdes}
   \Des(u)\ \subseteq\ \setS \ \iff\
    u\ \leq\ Z(\setS)\ =\ \zeta_\setS\,.
 \end{equation}
The notion of global descents is a very natural companion of that of
(ordinary) descents, in that the map
$\GDes\colon\frakS_n\to\calQ_n$
is {\em right} adjoint to $Z\colon\calQ_n\to\frakS_n$.

\begin{prop}\label{P:galoisglobal}
 The pair of maps $(Z,\GDes):\calQ_n\rightleftarrows\frakS_n$
 is a Galois connection.
\end{prop}

\begin{proof}
 We already know that $Z$ is order preserving.
 So is $\GDes$, because $p$ is a global descent of a permutation $u$ if and
 only if $(i,j)\in\Inv(u)$ for every $i\leq p$, $j\geq p{+}1$.
 It remains to check that
 \begin{equation}\label{E:galoisgdes}
   \zeta_\setS\ \leq\  u\ \iff\ \setS\ \subseteq\ \GDes(u)\,.
 \end{equation}
As in the proof of Proposition~\ref{P:galois}, this follows from
 \begin{align*}
  \zeta_\setS\ &=\ \min\{u\in\frakS_n\mid \GDes(u)\subseteq\setS\}\\
   &=\ \min\{u\in\frakS_n\mid \GDes(u)=\setS\}\,,
 \end{align*}
 which is clear from the definition of $\zeta_\setS$.
\end{proof}

We turn to properties of the decomposition of $\frakS_n$ into left cosets
of $\frakS_\setS$ related to the notion of global descents.
Recall that $\Sh{(p,q)}$ is a set of representatives for the left cosets of
$\frakS_p\times\frakS_q$ in $\frakS_{p+q}$,
and that $\zeta_{p,q}=(q{+}1, q{+}2,\ldots, q{+}p,\, 1, 2, \ldots, q)$.
\begin{lemm}\label{L:global}
 Suppose $p,q$ are non-negative integers and let $w\in\frakS_{p+q}$.
 Then
 \[
   p\in\GDes(w) \iff
   w\equiv\zeta_{p,q} \mod \frakS_p\times\frakS_q\iff w\geq \zeta_{p,q}\,.
 \]
\end{lemm}
\begin{proof}
 First suppose that $w\in\frakS_{p+q}$ is in the same left coset of
 $\frakS_p\times\frakS_q$ as is $\zeta_{p,q}$.
 Thus, there are permutations $u\in\frakS_p$ and $v\in\frakS_q$ such that
 \[
    w\ =\ \zeta_{p,q}\cdot(u\times v)\,.
 \]
 If $i\leq p$, then $u_i\in\{1,\ldots,p\}$ and thus
 $w_i\in\{q+1,\ldots,q+p\}$, so $p$ is a global descent of $w$ as needed.

 For the other direction, suppose $p$ is a global descent of $w$ and set
 \[
   \overline{w}\ :=\ \zeta_{p,q}^{-1} \cdot w\ =\
    (p{+}1,p{+}2,\ldots,p{+}q\,,1,2,\ldots,p)\cdot w\,.
 \]
 Let $1\leq i\leq p$. 
 By assumption, $w_i\in\{q{+}1,\ldots,q{+}p\}$. 
 Hence $\overline{w}_i\in\{1,\ldots,p\}$, which means that 
 $\overline{w}=u\times v$ for some $u\in\frakS_p$ and $v\in\frakS_q$, as
 needed.
 
 Noting that $\zeta_{p,q}$ is a minimal coset representative and that the map
 $\lambda^{-1}$ is order preserving (Proposition~\ref{P:shuffles}(a)) proves
 the second equivalence.
\end{proof}

For any subset $\setS$ of $[n{-}1]$, we have the left coset map
$\lambda\colon\Sh{\setS}\times\frakS_\setS\to \frakS_n$ of
Section~\ref{S:multiweak}. 
Given a permutation $u\in\frakS_n$, consider its `projection' 
$u_\setS$ to $\frakS_\setS$, which is defined to be the second
component of $\lambda^{-1}(u)$.
That is, $\lambda^{-1}(u)=(\zeta,u_\setS)$ for some permutation 
$\zeta\in\Sh{\setS}$.
If $\setS=\{p_1<p_2<\cdots<p_k\}$, then  by Lemma \ref{L:stdecom},
 \begin{equation}\label{E:u^S}
   u_\setS\ =\ \st(u_1,\dotsc,u_{p_1})\times
               \st(u_{p_1+1},\dotsc, u_{p_2})\times\dots\times
               \st(u_{p_k+1},\dotsc, u_{n})\,.
 \end{equation}
In particular, $u_{\emptyset}=u$ and $u_{[n{-}1]}=1_n$. 
We relate this projection to the order and lattice structure of $\frakS_n$.
For $i<j$, let $[i,j):=\{i,i{+}1,\ldots,j{-}1\}$.

\begin{lemm}\label{L:invus}
  For any $u\in\frakS_n$ and subset $\setS$ of\/ $[n{-}1]$, 
 \[
   \Inv(u_\setS)\ =\ \{(i,j)\in\Inv(u)\mid [i,j)\cap \setS=\emptyset\}\,.
 \]
 In particular, $u_\setS\leq u$.
\end{lemm}

\begin{proof} 
 Let $i<j$ be integers in $[n]$. 
 Suppose that there is an element $p\in \setS$ with $i\leq p<j$. 
 Since $\frakS_\setS\subseteq\frakS_p\times\frakS_{n-p}$, we have 
 $u_\setS\in\frakS_p\times\frakS_{n-p}$, and so $u_\setS(i)<u_\setS(j)$.
 Thus $(i,j)\not\in\Inv(u_\setS)$. 
 Suppose now that that $[i,j)\cap \setS=\emptyset$.
 Then there are consecutive elements $p$ and $q$ of $\setS$ such that
 $p<i<j\leq q$. 
 By~\eqref{E:u^S}, 
 \[
   u_\setS(i)\ =\ p+\st(u_{p+1},\ldots,u_q)(i)\ \text{ and }\ 
   u_\setS(j)\ =\ p+\st(u_{p+1},\ldots,u_q)(j)\,.
 \]
 By \eqref{E:st}, this implies that 
 \[
   u_\setS(i)\ >\ u_\setS(j)\  \iff\  u(i)\ >\ u(j),
 \]
 and thus $(i,j)\in \Inv(u_\setS)\iff (i,j)\in\Inv(u)$.
 This completes the proof.
\end{proof}

\begin{prop}\label{P:latticeus}
 Let $u$, $v\in\frakS_n$ and\/ $\setS,\setT$ be subsets of\/ $[n{-}1]$.
 Then
 \begin{enumerate}
  \item[({\it i})]
        If $u\leq v$ then $u_\setS\leq v_\setS$ and if\/
        $\setT\subseteq\setS$ then $u_\setT\geq u_\setS$.
  \item[({\it ii})]
        $u_\setS\wedge v_\setT=(u\wedge v)_{\setS\cup\setT}$,
  \item[({\it iii})]
        If\/ $\setS\subseteq \GDes(v)$ and\/ $\setT\subseteq \GDes(u)$, then  
        $u_\setS\vee v_\setT=(u\vee v)_{\setS\cap\setT}$.
 \end{enumerate}
\end{prop}

\begin{proof} 
The first statement is an immediate consequence of Lemma~\ref{L:invus}.
For the second, we use~\eqref{E:defmeet} to show that  
$\Inv(u_\setS\wedge v_\setT)=\Inv((u\wedge v)_{\setS\cup\setT})$.

First, suppose $(i,j)\in\Inv((u\wedge v)_{\setS\cup\setT})$.
Then by Lemma~\ref{L:invus} and~\eqref{E:defmeet}, we have 
$[i,j)\cap \bigl(\setS\cup\setT\bigr)=\emptyset$, and given a chain 
$i=k_0<\dotsb<k_s=j$, there is an index $r$ such that
$(k_{r-1},k_r)\in\Inv(u)\cap\Inv(v)$. 
Hence we also have $[k_{r-1},k_r)\cap\bigl(\setS\cup\setT\bigr)=\emptyset$, and 
thus $(k_{r-1},k_r)\in\Inv(u_\setS)\cap\Inv(v_\setT)$. 
Thus $(i,j)\in\Inv(u_\setS\wedge v_\setT)$.

We show the other inclusion.
Let $(i,j)\in\Inv(u_\setS\wedge v_\setT)$. 
Considering the chain $i<j$, we must have 
$(i,j)\in\Inv(u_\setS)\cap\Inv(v_\setT)$.
In particular, $[i,j)\cap\bigl(\setS\cup\setT\bigr)=\emptyset$.
On the other hand, for any chain $i=k_0<\dotsb<k_s=j$ there is an index $r$ 
such that $(k_{r-1},k_r)\in\Inv(u_\setS)\cap\Inv(v_\setT)$.
Since this is a subset of $\Inv(u)\cap\Inv(v)$, we have 
$(i,j)\in\Inv(u\wedge v)$.
Together with $[i,j)\cap \setS\cup\setT=\emptyset$, we see that
$(i,j)\in\Inv((u\wedge v)_{\setS\cup\setT})$, proving the second
statement.\smallskip

For the third statement, first note that statement ({\it i}) implies that
$u_\setS\leq(u\vee v)_\setS\leq (u\vee v)_{\setS\cap\setT}$ 
and similarly $v_\setT\leq (u\vee v)_{\setS\cap\setT}$.
Thus we have $u_\setS\vee v_\setT \leq (u\vee v)_{\setS\cap\setT}$.
To show the other inequality, we need the assumptions on $\setS$ and $\setT$.
With those assumptions, we show $\Inv((u\vee v)_{\setS\cap\setT})
              \subseteq\Inv(u_\setS\vee v_\setT)$.

Suppose that $\setS\subseteq \GDes(v)$ and $\setT\subseteq \GDes(u)$, so that
$\setS$ consists of global descents of $v$ and 
$\setT$ consists of global descents of $u$.
Let $(i,j)\in\Inv((u\vee v)_{\setS\cap\setT})$. 
Then, by Lemma~\ref{L:invus} and~\eqref{E:defjoin},
$[i,j)\cap \setS\cap\setT=\emptyset$ and there is a chain 
$i=k_0<\dotsb<k_s=j$ such that for every $r$,
$(k_{r-1},k_r)\in\Inv(u)\cup\Inv(v)$. 
We refine this chain so that every pair of consecutive elements
belongs to $\Inv(u_\setS)\cup\Inv(v_\setT)$.

If $[k_{r-1},k_r)\cap (\setS\cup\setT)=\emptyset$ then, 
by Lemma~\ref{L:invus}, $(k_{r-1},k_r)\in\Inv(u_\setS)\cup\Inv(v_\setT)$ 
and this interval need not be refined. 
If however the intersection is not empty, then choose any refinement 
 \[
   k_{r-1}\ =\ k_0^{(r)}\ <\ k_1^{(r)}\ <\ \dotsb\ <\ k_{s_r}^{(r)}\ =\ k_r\,,
 \]
with the property that each interval $[k_{t-1}^{(r)},k_t^{(r)})$ contains 
exactly {\em one} element of $\setS$ or $\setT$, but not an element
of both. 
This is possible because $[i,j)\cap \setS\cap\setT=\emptyset$. 
We claim that each pair $(k_{t-1}^{(r)},k_t^{(r)})$ is in 
$\Inv(u_\setS)\cup\Inv(v_\setT)$. 
In fact, if $[k_{t-1}^{(r)},k_t^{(r)})$ contains an element $p\in \setS$, 
then that is a global descent of $v$, so 
$(k_{t-1}^{(r)},k_t^{(r)})\in\Inv(v)$.
Thus $(k_{t-1}^{(r)},k_t^{(r)})\in\Inv(v_\setT)$, since
$[k_{t-1}^{(r)},k_t^{(r)})\cap \setT=\emptyset$ by our construction 
of the refinement. 
Similarly, if $[k_{t-1}^{(r)},k_t^{(r)})$ contains an element of $\setT$,
 then $(k_{t-1}^{(r)},k_t^{(r)})\in\Inv(u_\setS)$.
We have thus constructed a chain from $i$ to $j$ with the required property,
which shows that $(i,j)\in\Inv(u_\setS\vee v_\setT)$ and completes the proof.
\end{proof}

We calculate the descents and global descents of some
particular permutations.
The straightforward proof is left to the reader.

\begin{lemm}\label{L:over/under}
Let $u\in\frakS_p$ and $v\in\frakS_q$. Then
\begin{align*}
({\text i}) & &\Des(u\times v) &\ =\ \Des(u)\cup\bigl(p+\Des(v)\bigr),\\
({\text ii}) & &\GDes(u\times v) &\ =\emptyset,\\ 
({\text iii}) & &\Des\bigl(\zeta_{p,q}\cdot(u\times v)\bigr) 
        &\ =\ \Des(u)\cup\{p\}\cup\bigl(p+\Des(v)\bigr),\\
({\text iv}) & &\GDes\bigl(\zeta_{p,q}\cdot(u\times v)\bigr) 
       &\ =\ \GDes(u)\cup\{p\}\cup\bigl(p+\GDes(v)\bigr)\hspace{30pt}
\end{align*}
More generally, let $u_{(i)}\in\frakS_{p_i}$, $i=1,\ldots,k$,
$\setS=\{p_1,p_1+p_2,\ldots,p_1+\cdots+p_{k-1}\}\subseteq[n{-}1]$. Then
\begin{align*}
({\text v}) & &\Des(u_{(1)}\times\cdots\times u_{(k)}) 
          &\ =\ \bigcup_{i=1}^{k} \bigl(p_1+\cdots+p_{i-1}+\Des(u_{(i)})\bigr),\\
({\text vi}) & &\GDes(u_{(1)}\times\cdots\times u_{(k)}) 
            &\ =\ \emptyset,\\ 
({\text vii}) & &
      \Des\bigl(\zeta_\setS\cdot(u_{(1)}\times\cdots\times u_{(k)})\bigr) 
      &\ =\ \setS\cup\bigcup_{i=1}^{k}
            \bigl(p_1+\cdots+p_{i-1}+\Des(u_{(i)})\bigr),\\
({\text viii}) & &
       \GDes\bigl(\zeta_\setS\cdot(u_{(1)}\times\cdots\times u_{(k)})\bigr) 
      &\ =\ \setS\cup\bigcup_{i=1}^{k}
            \bigl(p_1+\cdots+p_{i-1}+\GDes(u_{(i)})\bigr).
 \end{align*}
\end{lemm}

\begin{lemm}\label{L:higherglobal}
 Let $u\in\frakS_n$ and $\setS\subseteq [n{-}1]$.
 Then
 \[
   \setS \subseteq \GDes(u)\ \iff\  u\;=\;\zeta_\setS u_\setS\,.
 \]

\end{lemm}
\begin{proof}
The reverse implication follows from Lemma~\ref{L:over/under}~({\it viii}). The
other follows by induction from
Lemma~\ref{L:global} and~\eqref{E:maxshuffles}.
\end{proof}

\section{The coproduct of $\SSym$}\label{S:coproduct}

The coproduct of $\SSym$~\eqref{E:cop-malvenuto} takes a simple form on
the monomial basis.
We derive this formula using some results of
Section~\ref{S:bruhat}.
For a permutation $u\in\frakS_n$, define $\AGDes(u)$ to be
$\GDes(u)\cup\{0,n\}$.

\begin{thm}\label{T:cop-monomial}
 Let $u\in \frakS_n$.
 Then
 \begin{equation}\label{E:cop-monomial}
    \Delta(\calM_u)\ =\ \sum_{p\in\AGDes(u)}
    \calM_{\st(u_1,\ldots,u_p)}\ten\calM_{\st(u_{p+1},\ldots,u_n)}\,.
 \end{equation}
\end{thm}
\begin{proof}
 Let $\Delta'\colon \SSym\to\SSym\ten\SSym$ be the  map
 whose action on the monomial basis is defined by the
 sum~\eqref{E:cop-monomial}. 
 We show that $\Delta'$ is the coproduct $\Delta$, as defined
 by~\eqref{E:cop-malvenuto}. 
 We use the following notation.
 For $w\in\frakS_n$ and $0\leq p\leq n$, let 
 $w^p_{(1)}:=\st(w_1,\ldots,w_p)$ and 
 $w^p_{(2)}:=\st(w_{p+1},\ldots,w_n)$. 
 By virtue of Lemmas \ref{L:stdecom} and \ref{L:global}, we have
\[v\ =\ \zeta_{p,n-p}\cdot(v^p_{(1)}\times v^p_{(2)}) \iff p\in\AGDes(v)\,.\]
 Therefore,
 \begin{align*}
  \Delta'(\calF_u) & 
    = \sum_{u\leq v}\Delta'(\calM_v)=\sum_{u\leq v}\sum_{p\in\AGDes(v)}
          \calM_{v^p_{(1)}}\ten\calM_{v^p_{(2)}} \\
   & =\sum_{p=0}^n \sumsub{u\leq
       v\\v=\zeta_{p,n-p}\cdot(v^p_{(1)}\times v^p_{(2)})\rule{0pt}{10pt}}
           \calM_{v^p_{(1)}}\ten\calM_{v^p_{(2)}} 
     =\sum_{p=0}^n \sumsub{v_1,\,v_2
       \\u\leq\zeta_{p,n-p}\cdot(v_1\times v_2) \rule{0pt}{10pt}}
      \calM_{v_1}\ten\calM_{v_2}\,.
 \end{align*}
Write $u=\zeta\cdot(u^p_{(1)}\times u^p_{(2)})$ for some $\zeta\in\Sh{(p,n-p)}$
 which depends on $p$. By Proposition
\ref{P:shuffles},
\[ \zeta\cdot(u^p_{(1)}\times u^p_{(2)})\leq \zeta_{p,n-p}\cdot(v_1\times v_2) \iff
u^p_{(1)}\leq v_1 \text{ and }u^p_{(2)}\leq v_2\,.\]
Therefore,
 \begin{align*}
  \Delta'(\calF_u) & 
     =\sum_{p=0}^n \sumsub{v_1\,,v_2
       \\u^p_{(1)}\leq v_1\,,u^p_{(2)}\leq v_2}
           \calM_{v_1}\ten\calM_{v_2} 
     =\sum_{p=0}^n \sum_{u^p_{(1)}\leq v_1}\calM_{v_1} \otimes
\sum_{u^p_{(2)}\leq v_2}\calM_{v_2}\\
&=\sum_{p=0}^n
\calF_{u^p_{(1)}}\ten\calF_{u^p_{(2)}}=\Delta(\calF_u)\,.
 \end{align*}
\end{proof}

\begin{rem} \label{R:loday}
 The action of the coproduct of $\SSym$ on the fundamental basis can also be
 expressed in terms of the weak order.
 To see this, let $u\in\frakS_n$ and $0\leq p \leq n$ and write 
 $u=\zeta\cdot(u^p_{(1)}\times u^p_{(2)})$.
 By Proposition~\ref{P:shuffles},   $u^p_{(1)}\times u^p_{(2)}\leq u \leq 
           \zeta_{p,n-p}\cdot(u^p_{(1)}\times u^p_{(2)})$.
 Moreover, $u^p_{(1)}$ and $u^p_{(2)}$ are the only permutations in 
 $\frakS_p$ and $\frakS_{n-p}$ with this property, again by
 Proposition~\ref{P:shuffles}.  
 Therefore, equation~\eqref{E:cop-malvenuto} is also described by
 $\Delta(\calF_u)=\sum\calF_v\ten\calF_w$, where the sum is over all
 $p$ from 0 to $n$ and all permutations 
 $v\in\frakS_p$ and $w\in\frakS_{n-p}$ such that 
 $v\times w\leq u\leq \zeta_{p,n-p}\cdot(v\times w)$. 
 This fact (in its dual form) is due to Loday and
 Ronco~\cite[Theorem~4.1]{LR01}, who were the first to point out the relevance
 of the weak order to the Hopf algebra structure of $\SSym$. 
\end{rem}

\section{The product of $\SSym$}\label{S:product}
We give an explicit formula for the product of $\SSym$ in terms of its
monomial basis  and a geometric interpretation for the 
structure constants.
Remarkably, these are still non-negative integers.
For instance, 
 \begin{multline}\label{E:prod-M}\ 
  \calM_{12} \cdot\calM_{21}\ =\ \calM_{4312}+\calM_{4231}+\calM_{3421}+\calM_{4123}+\calM_{2341}\\
    +\calM_{1243}+\calM_{1423}+\calM_{1342}+3\calM_{1432}+
     2\calM_{2431}+2\calM_{4132}\,. \ 
 \end{multline}
The structure constants count special ways of shuffling two permutations,
according to certain conditions involving the weak order. 
Specifically, for $u\in\frakS_p$, $v\in\frakS_q$ and $w\in\frakS_{p+q}$,
define $A^w_{u,v}\subseteq\Sh{(p,q)}$ to be those $\zeta\in\Sh{(p,q)}$
satisfying
 \begin{equation}\label{E:A-def}
    \begin{array}{ll}
       (i) & (u\times v)\cdot\zeta^{-1} \leq w, \text{ and}\\
      (ii) & \text{if $u\leq u'$ and $v\leq v'$ satisfy
             $(u'\times v')\cdot\zeta^{-1}\leq w$,}\\
           & \text{then $u=u'$ and $v=v'$.}
     \end{array}
 \end{equation}
Set $\alpha^w_{u,v}:=\#A^w_{u,v}$.
We will prove the following theorem.

\begin{thm}\label{T:prod-monomial}
  For any $u\in\frakS_p$ and $v\in\frakS_q$, we have
 \begin{equation}\label{E:prod-monomial}
    \calM_u\cdot\calM_v\ =\ \sum_{w\in\frakS_{p+q}}\alpha^w_{u,v}\,\calM_w\,.
 \end{equation}
\end{thm}

For instance, in~\eqref{E:prod-M} the coefficient of $\calM_{2431}$ in 
$\calM_{12}\cdot\calM_{21}$ is 2 because among the six permutations
in $\Sh{(2,2)}$, 
 \[
    1234,\ 1324,\ 1423,\ 2314,\ 2413,\ 3412\,,
 \]
only the first two satisfy conditions $(i)$ and $(ii)$ of~\eqref{E:A-def}.
In fact, $2314$, $2413$ and $3412$ do not
satisfy $(i)$, while $1423$ satisfies $(i)$ but not $(ii)$.

The structure constants $\alpha^w_{u,v}$ admit a geometric-combinatorial
description in terms of the permutahedron. To derive it, recall the convex
embeddings of Proposition~\ref{P:rightshuffling}.
 \[
   \rho_\zeta\ \colon\ \frakS_p\times\frakS_q\ \to\ \frakS_{p+q}\,,
    \qquad \rho_\zeta(u,v)\ :=\ (u\times v)\cdot\zeta^{-1}\,.
 \]
Since $\rho_\zeta$ preserves joins, we may further rewrite the
definition~\eqref{E:A-def} of $A^w_{u,v}$  as
 \begin{equation}\label{E:defalpha}
   A^w_{u,v}\ =\ \bigl\{\zeta \in \Sh{(p,q)}\mid
                   (u,v)=\max\rho_\zeta^{-1}[1,w]\bigr\}\,,
 \end{equation}
where $[w,w']:=\{w''\mid w\leq w''\leq w'\}$ denotes
the interval between $w$ and $w'$.

\begin{figure}[htb]
 \[
   \setlength{\unitlength}{0.8pt}
   \begin{picture}(230,215)(-10,0)
    \put(0,0){\epsfysize=168.6pt\epsffile{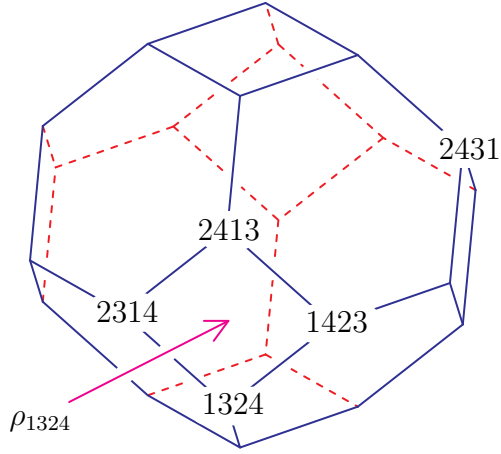}}
    \put(193,138){2431}
    \put(79,99.5){2413}
    \put(31,61){2314}  \put(130,56){1423}
    \put(81,17.5){1324}
    \put(-10,12){$\rho_{1324}$}
   \end{picture}
 \]
  \caption{The facet $\rho_{1324}$ of type $(2,2)$ and
    $w=2431$\label{F:alphaperm}.} \end{figure}

The vertices of the $(n{-}1)$-dimensional permutahedron can be   
indexed by the elements of $\frakS_n$ so that its 1-skeleton
is the Hasse diagram of the weak order (see Figure~\ref{F:S4}).
Facets of the permutahedron are products of two lower dimensional
permutahedra, and the image of $\rho_\zeta$ is the set of
vertices in a facet.
Moreover, every facet arises in this way for a unique triple 
$(p,q,\zeta)$ with $p+q=n$ and $\zeta\in \Sh{(p,q)}$;
see~\cite[Lemma~4.2]{Mi66}, or~\cite[Exer.~2.9]{BLSWZ},
or~\cite[Prop.~A.1]{Lod}. Let us say that such a 
facet has {\em type} $(p,q)$. Figure~\ref{F:alphaperm} displays the image of
$\rho_{1324}$, a facet of the $3$-permutahedron of type $(2,2)$, and the
permutation $2431$.

The description~\eqref{E:defalpha} of $A^w_{u,v}$ (and hence of $\alpha^w_{u,v}$) 
can be interpreted as follows: 
Given $u\in\frakS_p$, $v\in\frakS_q$, and
$w\in\frakS_{p+q}$, the structure constant $\alpha^w_{u,v}$
counts the number of facets of type $(p,q)$ of the $(p{+}q{-}1)$-permutahedron 
such that the vertex $\rho_\zeta(u,v)$ is below $w$ and it is the maximum
vertex in that facet below $w$.

For instance, the facet $\rho_{1324}$ contributes to the structure constant
$\alpha^{2431}_{12,21}$ because the vertex $\rho_{1324}(12,21)=1423$ satisfies
the required properties in relation to the 
vertex $w=2431$, as shown in Figure~\ref{F:alphaperm}.

This description of the product of $\SSym$ has an analog for $\QSym$ that we
present in Section~\ref{S:descentmap}.

\begin{proof}[Proof of Theorem~\ref{T:prod-monomial}]
 Expand the product $\calM_u\cdot\calM_v$ in the fundamental basis and then 
 use Formula~\eqref{E:prod-fundamental} to obtain
 \begin{align*}
   \calM_u\cdot\calM_v
    &\ =\  \sumsub{u\leq u'\\v\leq v'}
                \mu(u,u')\mu(v,v')\,\calF_{u'}\cdot\calF_{v'}\\
    &\ =\ \sum_{\zeta\in \Sh{(p,q)}}\sumsub{u\leq u'\\v\leq v'}
           \mu(u,u')\mu(v,v')\calF_{(u'\times v')\cdot\zeta^{-1}}\,.
 \end{align*}
 Expressing this result in terms of the monomial basis  and
 collecting like terms gives
 \begin{align*}
    \calM_u\cdot\calM_v
    &\ =\ \sum_{\zeta\in \Sh{(p,q)}}
       \sumsub{u\leq u',\ v\leq v'\\(u'\times v')\cdot\zeta^{-1}\leq w}
         \mu(u,u')\mu(v,v')\calM_{w}\\
     &\ =\ \sum_w
     \sumsub{u\leq u'\\v\leq v'}\mu(u,u')\mu(v,v')\beta^w_{u',v'}\calM_{w}\,,
 \end{align*}
 where $\beta^w_{u',v'}$ is the number of permutations in the set
 \[
    B^w_{u',v'}\ :=\ \bigl\{\zeta\in \Sh{(p,q)}\mid 
                    (u'\times v')\cdot\zeta^{-1}\leq w \bigr\}\,.
\]
The theorem will follow once we show that 
 \[
   \alpha^w_{u,v}\ =\ \sum_{u\leq u',\ v\leq v'}
         \mu(u,u')\mu(v,v')\beta^w_{u',v'}\,,
 \]
or equivalently, by M\"obius inversion on $\frakS_p\times\frakS_q$,
 \[
   \beta^w_{u,v}\ =\ \sum_{u\leq u',\ v\leq v'}\alpha^w_{u',v'}\,.
 \]
 We prove this last equality by showing that 
 \[
    B^w_{u,v}\ =\ \coprod_{u\leq u',\ v\leq v'}A^w_{u',v'}\;,
 \]
 where the union is disjoint.

 To see this, first suppose $\zeta\in A^w_{u,v}\cap A^w_{u',v'}$.  
 Then, by condition $(i)$ of~\eqref{E:A-def},
 \[
   (u\times v)\cdot\zeta^{-1}\leq w \quad\text{ and }\quad
   (u'\times v')\cdot\zeta^{-1}\leq w\,.
 \]
 By Proposition~\ref{P:rightshuffling}.(d), 
 \[
   \bigl((u\vee u')\times(v\vee v')\bigr)\cdot\zeta^{-1}\leq w\,.
 \]
 But then, by condition $(ii)$ of~\eqref{E:A-def},
 \[
    u\ =\ u\vee u'\ =\ u' \quad \text{ and }\quad 
    v\ =\ v\vee v'\ =\ v'\,,
 \]
 so the union is disjoint.

 Next, suppose that $\zeta\in A^w_{u',v'}$ for some $u\leq u'$ and 
 $v\leq v'$. 
 Then, by condition $(i)$ of~\eqref{E:A-def}, 
 $(u'\times v')\cdot\zeta^{-1}\leq w$. 
 By Proposition~\ref{P:rightshuffling}.(c) we have, 
 $(u\times v)\cdot\zeta^{-1}\leq w$, so $\zeta\in B^w_{u,v}$.
 This proves one inclusion.

 For the other inclusion, suppose that $\zeta\in B^w_{u,v}$. 
 Define
 \[
   (\overline{u},\overline{v})\ :=\ 
    \bigvee\{(u',v')\mid  u\leq u',\ v\leq v', \text{ and }
                          (u'\times v')\cdot\zeta^{-1}\leq w\}\,.
 \]
 Then $\zeta\in A^w_{\overline{u},\overline{v}}$: condition $(i)$ is satisfied
 because $\rho_\zeta$ preserves joins, and $(ii)$ simply by definition.
 This completes the proof.
\end{proof}

\section{The antipode of $\SSym$}\label{S:antipode}
Malvenuto left open the problem of an explicit formula for the antipode of
$\SSym$~\cite[pp.~59--60]{Malv}. 
We solve that problem, giving formulas that identify the coefficients
of the antipode in terms of both the fundamental and monomial basis in
explicit combinatorial terms. 

We first review a general formula for the antipode of a connected Hopf 
algebra $H$, due to Takeuchi~\cite[Lemma 14]{Tak71}
(see also Milnor and Moore~\cite{MM65}). 
Let $H$ be an arbitrary bialgebra  with structure maps:
multiplication $m\colon H\ten H\to H$, 
unit $u\colon \Q\to H$, 
comultiplication $\Delta\colon H\to H\ten H$, and
counit $\epsilon\colon H\to \Q$. 
Set $m^{(1)}=m$, $\Delta^{(1)}=\Delta$, and for any $k\geq 2$,
 \begin{align*}
   m^{(k)}    &\ =\ m(m^{(k-1)}\ten \id)\ \colon\
                     H^{\ten k+1}\to H,\,\quad \text{ and}\\
  \Delta^{(k)}&\ =\ (\Delta^{(k-1)}\ten \id)\Delta\ \colon\ H\to H^{\ten k+1}\,.
 \end{align*}
These are the higher or iterated products and coproducts. 
We also set
 \begin{align*}
  m^{(-1)}     &\ =\ u\colon \Q\to H,\\
  \Delta^{(-1)}&\ =\ \epsilon\colon H\to \Q,\quad  \text{ and }\\
  m^{(0)}      &\ =\ \Delta^{(0)}\ =\ \id\colon H\to H\,.
 \end{align*}
If $f:H\to H$ is any linear map, the convolution powers of $f$ are, for any $k\geq 0$,
 \[
   f^{\ast k}\ =\ m^{(k-1)}f^{\ten k}\Delta^{(k-1)}\,.
 \]
In particular, $f^{\ast 0}=u\epsilon$ and $f^{\ast 1}=f$.

We set $\pi:=id-u\epsilon$. 
If $\pi$ is locally nilpotent with respect to convolution, then
$id=u\epsilon+\pi$ is invertible with respect to convolution, with inverse
\begin{equation}\label{E:antipode}
  S\ =\ \sum_{k\geq 0}(-\pi)^{\ast k}\ =\ \sum_{k\geq 0}(-1)^k
  m^{(k-1)}\pi^{\ten k}\Delta^{(k-1)}\,.
\end{equation}
This is certainly the case if $H$ is a graded connected bialgebra, 
in which case $\pi$ annihilates the component of degree $0$ 
(and hence $\pi^{\ast k}$ annihilates components of degree $<k$). 
Thus~\eqref{E:antipode} is a general formula for the antipode of
a graded connected Hopf algebra.

We will make use of this formula to find explicit formulas for the antipode
of $\SSym$. 
The first task is to describe the higher products and coproducts
explicitly.
We begin with the higher coproducts in terms of the fundamental and
monomial bases.

\begin{prop}\label{P:highercop} 
 Let $v\in\frakS_n$, $n\geq 0$, and $k\geq 1$. Then
 \begin{itemize}  
  \item[({\it i})] \qquad\hspace{5.7pt}${\displaystyle 
        \Delta^{(k)}(\calF_v)\ =\ \sum_{0\leq p_1\leq\cdots\leq p_k\leq n}
                \calF_{\st(v_1,\,\ldots,\,v_{p_1})}\ten\dotsm\ten
               \calF_{\st(v_{p_k+1},\,\ldots,\,v_n)}}$, \ and 
   \item[({\it ii})]\qquad \rule{0pt}{17pt}${\displaystyle 
        \Delta^{(k)}(\calM_v)\  =\ 
    \sumsub{0\leq p_1\leq\ldots\leq p_k\leq n\\p_i\in\AGDes(v)\rule{0pt}{10pt}}
          \calM_{\st(v_1,\,\ldots,\,v_{p_1})}\ten\dotsm\ten
                 \calM_{\st(v_{p_k+1},\,\ldots,\,v_n)}}$.
 \end{itemize}
\end{prop}

\begin{proof} 
 Both formulas follow by induction from the corresponding descriptions of
 the coproduct, equations~\eqref{E:cop-malvenuto} and~\eqref{E:cop-monomial}.
\end{proof}

We describe higher products in terms of minimal coset representatives 
$\Sh{\setS}$ of parabolic subgroups, whose basic properties were discussed in
Section~\ref{S:multiweak}. 
Recall that for a subset $\setS=\{p_1<p_2<\cdots<p_k\}$ of $[n{-}1]$, we have
$\Sh{\setS}=\{\zeta\in\frakS_n\mid \Des(\zeta)\subseteq\setS\}$. 
Analogously to~\eqref{E:A-def}, given permutations $v_{(1)}\in\frakS_{p_1}$, 
$v_{(2)}\in\frakS_{p_2-p_1}$, \ldots, $ v_{(k{+}1)}\in\frakS_{n-p_k}$, 
define $A^w_{v_{(1)},v_{(2)},\ldots,v_{(k{+}1)}}\subseteq\Sh{\setS}$ 
to be those $\zeta\in\Sh{\setS}$ satisfying
 \begin{equation}\label{E:AS-def}
    \begin{array}{ll}
       (i) & \bigl(v_{(1)}\times v_{(2)}\times\dotsm\times
               v_{(k{+}1)}\bigr)\cdot\zeta^{-1} \leq w,\   \text{ and}\\
      (ii) & \text{if $v_{(i)}\leq v'_{(i)}\ \forall i$ and }\rule{0pt}{13pt}
             \bigl(v'_{(1)}\times v'_{(2)}\times\dotsm\times
               v'_{(k{+}1)}\bigr)\cdot\zeta^{-1} \leq w,\ \\
           & \text{then $v_{(i)}= v'_{(i)},\ \forall i$\,.}
     \end{array}
 \end{equation}
Set $\alpha^w_{v_{(1)},v_{(2)},\ldots,v_{(k{+}1)}}:=
     \#A^w_{v_{(1)},v_{(2)},\ldots,v_{(k{+}1)}}$.

\begin{prop}\label{P:higherprod} 
 Let $\setS$ and $v_{(1)},\ldots,v_{(k{+}1)}$ be as in the preceding paragraph.
 Then
 \begin{itemize}
    \item[({\it i})] 
           \qquad\hspace{13.2pt} ${\displaystyle
       \calF_{v_{(1)}}\cdot\calF_{v_{(2)}}\dotsm\calF_{v_{(k{+}1)}}\
       =\ \sum_{\zeta\in \Sh{\setS}} 
       \calF_{(v_{(1)}\times v_{(2)}\times\dotsm\times v_{(k{+}1)})
          \cdot\zeta^{-1}}}$ \ and 
    \item[({\it ii})]
         \qquad\rule{0pt}{17pt}${\displaystyle 
    \calM_{v_{(1)}}\cdot\calM_{v_{(2)}}\dotsm\calM_{v_{(k{+}1)}}\
     =\    \sum_{w\in\frakS_n}\alpha^w_{v_{(1)},v_{(2)},\ldots,v_{(k{+}1)}}\,
     \calM_w\,}$, 
 \end{itemize}
\end{prop}

\begin{proof} 
 The first formula follows immediately by
 induction  from \eqref{E:prod-fundamental} (the case $k=2$), 
 using~\eqref{E:highershuffles}.
 The second formula can be deduced
 from~({\it i}) in the same way as in the proof of
 Theorem~\ref{T:prod-monomial}.  
\end{proof}

The structure constants for the iterated product 
admit a geometric description similar to that of the product.
The image of the map
 \[
    \rho_\zeta\ \colon\ \frakS_\setS\to\frakS_n,\qquad  
     (v_{(1)}\times \dotsm\times v_{(k{+}1)})\ \longmapsto\  
     (v_{(1)}\times \dotsm\times v_{(k{+}1)}) \cdot\zeta^{-1}\,,
 \]
consists of the vertices of a face of codimension $k$ in the
$(n{-}1)$-permutahedron, and every such face arises in this way for a
unique pair $(\setS,\zeta)$ with $\setS\subseteq[n{-}1]$ having $k$ elements 
and $\zeta\in \Sh{\setS}$. 
Let us say that such a face has \emph{type} $\setS$.
The structure constant
$\alpha^w_{v_{(1)},\ldots,v_{(k{+}1)}}$ counts the number of faces of
type $\setS$  with the property that the vertex 
$\rho_\zeta(v_{(1)},\ldots,v_{(k{+}1)})$ is  below $w$ and it is
the maximum vertex in its face  below $w$.
\smallskip

We next determine the convolution powers of the projection $\pi=id-u\epsilon$. 
Recall that for any subset $\setS=\{p_1<p_2<\dotsb<p_k\}\subseteq[n{-}1]$ and 
$v\in\frakS_n$ we have
 \[
   v_\setS\ :=\ 
      \st(v_1,\ldots,v_{p_1})\times\st(v_{p_1+1},\ldots,v_{p_2})\times\dotsb\times
       \st(v_{p_k+1},\ldots,v_{n})\ \in\ \frakS_n\,,
\]
as given by~\eqref{E:u^S}.
We slightly amend our notation in order to simplify some
subsequent statements.
For $v,w\in\frakS_n$ and $\setS\subseteq[n{-}1]$, set
$A_\setS(v,w):=A^w_{v_{(1)}, \ldots, v_{(k+1)}}$, where
$v_{(1)},\ldots,v_{(k+1)}$ are the factors of $v_\setS$ in the definition
above. 
Comparing with~\eqref{E:AS-def}, we see that $A_\setS(v,w)\subseteq\Sh{\setS}$
consists of those $\zeta\in\Sh{\setS}$ satisfying 
 \begin{equation}\label{E:AwsetS}
   \begin{array}{ll}
      (i) & v_\setS\zeta^{-1}\leq w,\   \text{ and}\\ \rule{0pt}{13pt}
     (ii) & \text{if $v\leq v'$ and $v'_\setS\zeta^{-1}\leq w$ then $v=v'$} 
   \end{array}
 \end{equation}
Similarly, we define $\alpha_\setS(v,w):=\#A_\setS(v,w)$.
If $v_{(1)},\ldots,v_{(k+1)}$ are the factors in the definition of
$v_\setS$, then 
\begin{equation}\label{E:eq-consts}
  \alpha_\setS(v,w)\ =\ \alpha^w_{v_{(1)}, \ldots, v_{(k+1)}}\,.
\end{equation}

Let $\binom{[n{-}1]}{k{-}1}$ be the collection of subsets of $[n{-}1]$ of
size $k{-}1$.

\begin{prop}\label{P:pipower} 
 Let $n,k\geq 1$ and $v\in\frakS_n$. 
 Then
 \begin{itemize}
  \item[({\it i})] 
    \qquad\hspace{5.8pt}${\displaystyle
   \pi^{\ast k}(\calF_v)\ =\ \sum_{w\in\frakS_n}\ 
     \sumsub{\setS\in\binom{[n{-}1]}{k-1}\\\Des(w^{-1}v_\setS)\subseteq\setS}
           \calF_w}$, \ and 
    \item[({\it ii})] \rule{0pt}{18pt}\qquad
    ${\displaystyle
           \pi^{\ast k}(\calM_v)\ =\  \sum_{w\in\frakS_n}\ 
           \sumsub{\setS\in\binom{\GDes(v)}{k-1}}\alpha_\setS(v,w)\,\calM_w}$.
 \end{itemize}
\end{prop}

\begin{proof} 
 By Proposition~\ref{P:highercop}({\it i}), 
 \[
   \Delta^{(k-1)}(\calF_v)\ =\ \sum_{0\leq p_1\leq\dotsb\leq p_{k-1}\leq n}
           \calF_{\st(v_1,\ldots,v_{p_1})}\ten
            \calF_{\st(v_{p_1+1},\ldots,v_{p_2})}\ten\dotsb\ten
             \calF_{\st(v_{p_{k-1}+1},\ldots,v_{n})}\,.
 \]
 Suppose that an equality  $p_i=p_{i+1}$ occurs (where we define $p_0=0$ and
 $p_k=n$). 
 The corresponding permutation $\st(v_{p_i+1},\dotsc,v_{p_{i+1}})$ is then
 simply the unique permutation in $\frakS_0$, which indexes the element
 $1\in\ker(\pi)$.  
 Therefore, 
 \begin{align*}
   \pi^{\ast k}(\calF_v) 
   &\ =\ m^{(k-1)}\pi^{\ten k}\Delta^{(k-1)}(\calF_v)\\
   &\ =\ \sum_{0< p_1< p_2<\dotsb< p_{k-1}<n}
          \calF_{\st(v_1,\dotsc,v_{p_1})}\cdot
          \calF_{\st(v_{p_1+1},\dotsc,v_{p_2})}\dotsm
          \calF_{\st(v_{p_{k-1}+1},\dotsc,v_{n})}\\
   &\ =\ \sum_{0< p_1< p_2<\dotsb< p_{k-1}<n}\
         \sum_{\zeta\in\Sh{\{p_1,p_2,\dotsc,p_{k-1}\}}}
            \calF_{(\st(u_1,\dotsc,u_{p_1})\times\dotsb
                    \times\st(u_{p_{k-1}+1},\dotsc,u_n))\cdot\zeta^{-1}}\,,
\end{align*}
 the last equality by the formula of Proposition~\ref{P:higherprod}({\it i})
 for the iterated product. 
 Changing the index of summation in the first sum
 to $\setS\in\binom{[n{-}1]}{k-1}$ and using the definition of $v_\setS$ gives
 \[
    \pi^{\ast k}(\calF_v) 
   \ =\ \sum_{\setS\in\binom{[n{-}1]}{k-1}}
             \sum_{\zeta\in \Sh{\setS}} \calF_{v_\setS\zeta^{-1}}\,.
 \]
 Again reindexing the sum and using that $\Sh{\setS}$ consists of
 permutations whose descent set is a subset of $\setS$, we obtain
 \[
 \pi^{\ast k}(\calF_v) \ =\  
    \sum_{w\in\frakS_n}\sumsub{\setS\in\binom{[n{-}1]}{k-1}\\ 
            w^{-1}v_\setS\in \Sh{\setS}} \calF_w
   \ =\ \sum_{w\in\frakS_n}\sumsub{\setS\in\binom{[n{-}1]}{k-1}\\ 
                 \Des(w^{-1}v_\setS) \subseteq\setS}\calF_w\,,
 \]
 establishing ({\it i}).

 The second formula in terms of the monomial basis 
 follows in exactly the same manner from
 Propositions~\ref{P:highercop}({\it ii}) and~\ref{P:higherprod}({\it ii}) for
 the higher coproducts and products in terms of the monomial basis,
 using~\eqref{E:eq-consts}.
\end{proof}

We  derive explicit formulas for the antipode on both bases.
The formula for the fundamental basis is immediate from
Proposition~\ref{P:pipower}({\it i}) and~\eqref{E:antipode}.

\begin{thm}\label{T:ant-fundamental}
 For $v,w\in\frakS_n$ set
 \begin{align*}
     \lambda(v,w)\ :=\ &\ \#\{\setS\subseteq[n{-}1]\mid
     \Des(w^{-1}v_\setS)\subseteq\setS \text{ and $\#\setS$ is odd}\}\\
       & -\#\{\setS\subseteq[n{-}1]\mid \Des(w^{-1}v_\setS)\subseteq\setS 
         \text{ and $\#\setS$ is even}\}.
 \end{align*}
 Then
 \begin{equation}\label{E:ant-fundamental}
   S(\calF_v)\ =\ \sum_{w\in\frakS_n}\lambda(v,w)\,\calF_w\,.
 \end{equation}
\end{thm} 

The coefficients of the antipode on the fundamental basis may indeed be
positive or negative. For instance 
  \[
     S(\calF_{231})\ =\ \calF_{132}-\calF_{213}-2\calF_{231}+\calF_{312}\,.
 \]
The coefficient of $\calF_{312}$ is $1$ because $\{1\}$, $\{2\}$, and
$\{1,2\}$ are the subsets $\setS$ of $\{1,2\}$ which satisfy 
$\Des\bigl((312)^{-1}(231)_\setS\bigr)\subseteq\setS$.

Our description of these coefficients is semi-combinatorial, in the sense
that it involves a difference of cardinalities of sets. 
On the monomial basis the situation is different. 
The sign of the coefficients of $S(\calM_v)$ only depends on the number of
global descents of $v$. 
We provide a fully combinatorial description of these coefficients.
Let $v,w\in\frakS_n$ and suppose $\setS\subseteq\GDes(v)$.
Define $C_\setS(v,w)\subseteq\Sh{\setS}$ to be those $\zeta\in\Sh{\setS}$
satisfying 
 \begin{equation}\label{E:C-def}
   \begin{array}{rl}
       (i) & \text{$v_\setS\zeta^{-1}\leq w$,}\\
      (ii) & \text{if $v\leq v'$ and $v'_\setS\zeta^{-1}\leq w$ then $v=v'$, \
                 and}\\
     (iii) &  \text{if $\Des(\zeta)\subseteq\setR\subseteq\setS$ and 
              $v_\setR\zeta^{-1}\leq w$ then  $\setR=\setS$.}
   \end{array}
 \end{equation}
Set $\kappa(v,w):=\#C_{\GDes(v)}(v,w)$.

\begin{thm}\label{T:ant-monomial} 
  For $v,w\in\frakS_n$, we have
 \begin{equation}\label{E:ant-monomial}
   S(\calM_v)\ =\ (-1)^{\#\GDes(v)+1}\sum_{w\in\frakS_n}\kappa(v,w)\,\calM_w\,.
 \end{equation}
\end{thm}

For instance,
 \begin{multline*}
   S(\calM_{3412})\ =\ \calM_{1234}+2\calM_{1324}+\calM_{1342}+\calM_{1423}\\
    +\calM_{2314}+\calM_{2413}+\calM_{3124}+\calM_{3142}+\calM_{3412}\,.
 \end{multline*}
Consider the coefficient of $\calM_{3412}$. 
In this case, $\setS=\GDes(3412)=\{2\}$, so
 \[
   \Sh{\setS}\ =\ \{1234,\ 1324,\ 1423,\ 2314,\ 2413,\ 3412\}\,.
 \]
Then $1234$ satisfies ($i$) and ($ii$) of~\eqref{E:C-def} but not ($iii$),
1324 satisfies ($i$) and ($iii$) but not ($ii$), 
1423, 2314 and 2413 do not satisfy ($i$), and 3412 is the only 
element of $\Sh{\{2\}}$ that satisfies all three conditions of~\eqref{E:C-def}. 
Therefore $C_{\setS}(3412,3412)=\{3412\}$ and the coefficient
is $\kappa(3412,3412)=1$.

\begin{rem}
 The antipode of $\SSym$ has infinite order. 
 In fact, one may verify by induction that 
 \[
   S^{2m}(\calM_{231})\ =\ \calM_{231}+2m(\calM_{213}-\calM_{132})
    \quad\forall\ m\in\Z\,.
 \]
\end{rem}

\begin{proof}[Proof of Theorem~\ref{T:ant-monomial}]
 By formula~\eqref{E:antipode} and Proposition~\ref{P:pipower}({\it ii}), we have 
 \[
   S(\calM_v)\ =\ \sum_{w\in\frakS_n}\ \sum_{\setS\subseteq\GDes(v)}
             (-1)^{\#\setS+1}\alpha_\setS(v,w)\,\calM_w\,.
 \]
 For any $\setT\subseteq\GDes(v)$, define
 \begin{equation}\label{E:gamma_T}
   \gamma_\setT(v,w)\ :=\ \sum_{\setS\subseteq\setT}
       (-1)^{\#\setT\setminus\setS}\alpha_\setS(v,w)\ =\ 
       \sum_{\setS\subseteq\setT} \mu(\setS,\setT)\alpha_\setS(v,w)\,,
 \end{equation}
 where $\mu(\cdot,\cdot)$ is the M\"obius function of the Boolean poset
 $\calQ_n$. 
 We then have
 \[
   S(\calM_v)\ =\ (-1)^{\#\GDes(u)+1}\sum_{w\in\frakS_n}
                   \gamma_{\GDes(v)}(v,w)\calM_w\,.
 \]
 We complete the proof by showing that 
 $\kappa(v,w)=\gamma_{\GDes(v)}(v,w)$, and more generally that 
 $\gamma_\setS(v,w)=\#C_\setS(v,w)$, where $C_\setS(v,w)$ is defined
 in~\eqref{E:C-def}. 

 M\"obius inversion using the definition~\eqref{E:gamma_T} of 
 $\gamma_\setT(v,w)$ gives 
 \[
   \alpha_\setT(v,w)\ =\ \sum_{\setS\subseteq\setT}\gamma_\setS(v,w)\,.
 \]
 We prove this last equality by showing that
 \begin{equation}\label{E:disjoint-dec}
   A_\setT(v,w)\ =\ \coprod_{\setS\subseteq\setT} C_\setS(v,w)\,,
 \end{equation}
 where the union is disjoint.
 This implies that $\gamma_\setS(v,w)=\# C_\setS(v,w)$, which 
 will complete the proof.
 We argue that this is a disjoint union in several steps.\smallskip

 \noindent{\bf Claim 1:} If $\setS\subseteq\setT\subseteq\GDes(v)$ 
 then $A_\setS(v,w)\subseteq A_\setT(v,w)$.

 Let $\zeta\in A_\setS(v,w)$. 
 First of all, $\zeta\in \Sh{\setS}\subseteq \Sh{\setT}$, as $\Sh{\setS}$ is
 the set of permutations with descent set a subset of $\setS$.
 By condition ($i$) of~\eqref{E:C-def}, $v_\setS\zeta^{-1}\leq w$. 
 On the other hand, Proposition~\ref{P:latticeus}({\it i}) implies that 
 $u_\setT\leq u_\setS$ and
 both permutations are elements of the parabolic subgroup $\frakS_\setS$.
 Hence by Proposition~\ref{P:multirightshuffling}, 
 $u_\setT\zeta^{-1}\leq u_\setS\zeta^{-1}$. 
 Thus $u_\setT\zeta^{-1}\leq w$, which establishes condition ($i$)
 of~\eqref{E:C-def} for $\zeta$ to be in $A_\setT(v,w)$.

 Now suppose that $v\leq v'$ with $v'_\setT\zeta^{-1}\leq w$. 
 Since $v_\setS\zeta^{-1}\leq w$, we deduce that
 \[
    w\ \geq\ (v_\setS\zeta^{-1})\vee(v'_\setT\zeta^{-1})\ =\ 
      (v_\setS\vee v'_\setT)\zeta^{-1}\ =\ 
      (v\vee v')_{\setS\cap\setT}\zeta^{-1}\ =\ v_\setS\zeta^{-1}\,.
 \]
 The first equality is because $\rho_\zeta$ is a convex embedding and hence
 preserves joins by Proposition~\ref{P:multirightshuffling},
 and the second follows from Proposition~\ref{P:latticeus}({\it iii}) as 
 $\setS,\setT\subseteq\GDes(v)\subseteq\GDes(v')$.
 Hence, by condition ($ii$) for $A_\setS(v,w)$, we have $v=v'$. 
 This establishes ($ii$) for $\zeta$ to be in $A_\setT(v,w)$
 and completes the proof of Claim 1.\smallskip

 \noindent{\bf Claim 2:}
 If  $\setS, \setT \subseteq\GDes(v)$, then
 $A_{\setS}(v,w)\cap A_{\setT}(v,w)\ =\ A_{\setS\cap\setT}(v,w)$.

 The inclusion 
 $A_{\setS\cap\setT}(v,w)\subseteq A_{\setS}(v,w)\cap A_{\setT}(v,w)$ is a
 consequence of Claim 1. 
 To prove the converse, let 
 $\zeta\in A_{\setS}(v,w)\cap A_{\setT}(v,w)$. 
 Note that $\zeta\in \Sh{\setS}\cap \Sh{\setT}$, which equals
 $\Sh{\setS\cap\setT}$.

 By condition ($i$) for $\zeta\in A_{\setS}(v,w)$ and for 
 $\zeta\in A_{\setT}(v,w)$, we have $v_{\setS}\zeta^{-1}\leq w$ and 
 $v_{\setT}\zeta^{-1}\leq w$.
 Therefore,
 \[
   w\ \geq\ (v_{\setS}\zeta^{-1})\vee(v_{\setT}\zeta^{-1})\ =\ 
   (v_{\setS}\vee v_{\setT})\zeta^{-1}\ =\ 
   v_{\setS\cap\setT}\zeta^{-1}\,.
 \] 
 As before, this uses Proposition~\ref{P:latticeus}({\it iii}), which applies as 
 $\setS, \setT\subseteq\GDes(u)$. 
 This proves condition ($i$)
 of~\eqref{E:C-def} for $\zeta$ to be in $A_{\setS\cap\setT}(v,w)$.

 Now suppose that $v\leq v'$ with $v'_{\setS\cap\setT}\zeta^{-1}\leq w$. 
 By Proposition~\ref{P:latticeus}({\it i}), 
 $v'_{\setS}\leq v'_{\setS\cap\setT}$. 
 Then by Proposition~\ref{P:multirightshuffling},
 $v'_\setS\zeta^{-1}\leq v'_{\setS\cap\setT}\zeta^{-1}$. 
 Thus $v'_{\setS}\zeta^{-1}\leq w$ and by
 condition ($ii$) for $A_{\setS}(v,w)$ we deduce that $v=v'$. 
 This proves condition ($ii$) for  $\zeta$ to be in $A_{\setS\cap\setT}(v,w)$,
 and establishes Claim 2.\smallskip

 We complete the proof by showing that for $\setT\subseteq\GDes(v)$ we have the
 decomposition~\eqref{E:disjoint-dec} of $A_\setT(v,w)$ into disjoint subsets
 $C_\setS(v,w)$.
 Comparing the definitions~\eqref{E:AwsetS} and~\eqref{E:C-def}, we see that
 $C_\setS(v,w)\subseteq A_\setS(v,w)$. 
 Together with Claim 1 this implies that the right hand side
 of~\eqref{E:disjoint-dec} is contained in the left hand side.

 We show the union is disjoint.
 Suppose there is a permutation $\zeta\in C_{\setS}(v,w)\cap C_{\setS'}(v,w)$.  
 Then $\zeta\in A_{\setS}(v,w)\cap A_{\setS'}(v,w)$ which equals
 $A_{\setS\cap\setS'}(v,w)$, by Claim 2. 
 Hence, by condition ($i$) for $\zeta$ to be in $A_{\setS\cap\setS'}(v,w)$, we
 have $v_{\setS\cap\setS'}\zeta^{-1}\leq w$. 
 But then, from condition ($iii$) for $C_{\setS}(v,w)$ and for 
 $C_{\setS'}(v,w)$, we deduce that $\setS=\setS\cap\setS'=\setS'$, proving
 the union is disjoint.  

 We show that $A_\setT(v,w)$ is contained in the union
 of~\eqref{E:disjoint-dec}. 
 Let $\zeta\in A_\setT(v,w)$ and set
 \begin{equation}\label{E:setS}
   \setS\ :=\ \bigcap \{\setR\mid \setR\subseteq\setT,\ 
                        \zeta\in A_\setR(v,w)\}\,.
 \end{equation}
 By Claim 2, 
 \[
   A_\setS(v,w)\ =\ \bigcap \{A_\setR(v,w)\mid \setR\subseteq\setT,\ 
                              \zeta\in A_\setR(v,w)\}\,,
 \]
 so $\zeta\in A_\setS(v,w)$. 
 To show that $\zeta\in C_\setS(v,w)$, we must verify condition
 ($iii$) of~\eqref{E:C-def}.

 Suppose $\Des(\zeta)\subseteq\setR\subseteq\setS$ and 
 $v_\setR\zeta^{-1}\leq w$. 
 We need to show that $\setS\subseteq\setR$. 
 By the definition~\eqref{E:setS} of $\setS$, it suffices to show that
 $\zeta\in A_\setR(v,w)$. 
 By our assumption that $v_\setR\zeta^{-1}\leq w$, condition ($i$) for $\zeta$
 to be in $A_\setR(v,w)$ holds.  
 We show that condition ($ii$) also holds.
 Suppose $v\leq v'$ and $v'_\setR\zeta^{-1}\leq w$. 
 By Proposition~\ref{P:latticeus}({\it i}) we have
 $v'_{\setS}\leq v'_{\setR}$, and so by 
 Proposition~\ref{P:multirightshuffling}, 
 $v'_{\setS}\zeta^{-1}\leq v'_{\setR}\zeta^{-1}$. 
 Thus $v'_{\setS}\zeta^{-1}\leq w$, and by condition ($ii$)
 for $\zeta$ to be in $A_\setS(v,w)$, we have $v=v'$. 
 This establishes condition ($ii$) for $\zeta$ to be in $A_\setR(v,w)$. 
 Thus, $\zeta\in A_\setR(u,w)$, and as explained above, shows 
 that~\eqref{E:disjoint-dec} is a disjoint union and completes the 
 proof of the theorem.
\end{proof}

\section{Cofreeness, primitive elements, and the coradical filtration of $\SSym$}
\label{S:cofree}

The monomial basis reveals the existence of a second coalgebra grading on
$\SSym$,  given by the number of global descents of the indexing permutations.
We show that with respect to this grading, $\SSym$ is a cofree graded
coalgebra. We deduce an elegant description of the coradical
filtration: it corresponds to a filtration of the symmetric groups by certain
lower order ideals determined by the number of global descents.
In particular, the space of primitive elements is spanned by those 
$\calM_u$ where $u$ has no global descents.

We review the notion of cofree graded coalgebras. 
Let $V$ be a vector space and set
 \[ 
   Q(V)\ :=\ \bigoplus_{k\geq 0} V^{\ten k}\,.
 \]
The space $Q(V)$, graded by $k$, becomes a graded coalgebra with the
\emph{deconcatenation} coproduct
 \[
   \Delta(v_1\ten\ldots\ten v_k)\ =\ 
   \sum_{i=0}^k (v_1\ten\dotsb\ten v_i)\otimes(v_{i+1}\ten\dotsb\ten v_k)\,,
 \]
and counit $\epsilon(v_1\ten\dotsb\ten v_k)=0$ for $k\geq 1$.  $Q(V)$ is
connected, in the sense that the component of degree $0$ is identified with the
base field via $\epsilon$.

We call $Q(V)$ the cofree graded coalgebra cogenerated by $V$.
The canonical projection $\pi:Q(V)\to V$ satisfies the following universal
property. 
Given a graded coalgebra $C=\oplus_{k\geq 0}C^k$ and a linear map  
$\varphi:C\to V$ where $\varphi(C^k)=0$ when $k\neq 1$, 
there is a unique morphism of graded coalgebras
$\hat{\varphi}:C\to Q(V)$ such that the following diagram commutes
 \[
   \xymatrix{{\ C\ }\ar@{-->}[rr]^{\hat{\varphi}}\ar[dr]_{\varphi} &
   &{Q(V)}\ar[ld]^{\pi}\\ & {V} }
 \]
Explicitly, $\hat{\varphi}$  is defined by
 \begin{equation} \label{E:cofree}
   \hat{\varphi}_{|_{C^k}}=\varphi^{\ten k}\Delta^{(k-1)}\,.
 \end{equation}
In particular, $\hat{\varphi}_{|_{C^0}}=\epsilon$,
$\hat{\varphi}_{|_{C^1}}=\varphi$, and
$\hat{\varphi}_{|_{C^2}}=(\varphi\ten\varphi)\Delta$.

We establish the cofreeness of $\SSym$ by first defining
a second coalgebra grading.
Let $\frakS^{0}:=\frakS_0$, and for
$k\geq 1$, let
 \begin{align*}
   \frakS_n^{k}&\ :=\ \{u\in\frakS_n\mid u
                \text{ has exactly $k{-}1$ global descents}\}, \
       \text{ and }\\
    \frakS^{k}&\ :=\ \coprod_{n\geq 0}\frakS_n^{k}\,.
 \end{align*}
For instance,
 \begin{multline*}
  \frakS^{1}\ =\ \{1\}\ \cup\ \{12\}\ \cup\ \{123,\,213,\,132\}\
                 \cup\ \{1234,\,2134,\,1324,\,1243,\,3124,\\
                   2314,\,2143,\,1423,\,1342,\,3214,\,
                   3142,\,2413,\,1432\}\ \cup\ \dotsb
 \end{multline*}

Let $(\SSym)^k$ be the vector subspace of $\SSym$ spanned by
 $\{\calM_u\mid u\in\frakS^k\}$.

\begin{thm}\label{T:cofree} The decomposition $\SSym=\oplus_{k\geq 0}(\SSym)^k$
is a coalgebra grading. Moreover, endowed with this grading, $\SSym$ is a
cofree graded coalgebra. \end{thm}

\begin{proof}
Let $u\in\frakS_n^k$ and write $\GDes(u)=\{p_1<\ldots<p_{k-1}\}$. By
Theorem~\ref{T:cop-monomial},
 \[
    \Delta(\calM_u)\ =\ 1\ten\calM_u+\sum_{i=1}^{k-1}
    \calM_{\st(u_1,\ldots,u_{p_i})}\ten\calM_{\st(u_{p_i+1},\ldots,u_n)}+
    \calM_u\ten 1\, .
 \]
Since $\st(u_1,\ldots,u_{p_i})$ and $\st(u_{p_i+1},\ldots,u_n)$ have $i{-}1$
and $k{-}1{-}i$ global descents, we have
\[ \Delta\bigl((\SSym)^k\bigr)\ \subseteq\ 
 \bigoplus_{i=0}^k (\SSym)^i\ten (\SSym)^{k-i}\,.\]
Thus $\SSym=\oplus_{k\geq 0}(\SSym)^k$ is a graded coalgebra.

Let $V=(\SSym)^1$ and $\varphi:\SSym\to V$ the projection associated to the
grading. Let $\hat{\varphi}:\SSym\to Q(V)$ be the morphism of graded coalgebras
into the cofree graded coalgebra on $V$. For $u$ as above,
Proposition~\ref{P:highercop} gives,
 \[
   \Delta^{(k-1)}(\calM_u)\ =\
      \sumsub{0\leq q_1\leq \dotsb\leq q_{k-1}\leq n\\
              q_i\in\AGDes(u)\rule{0pt}{10pt}}
\calM_{\st(u_1,\dotsc,u_{q_1})}\ten\dotsm\ten
\calM_{\st(u_{q_{k-1}+1},\dotsc,u_{n})}\,.
 \]
Among these chains $0\leq q_1\leq \dotsb\leq q_{k-1}\leq n$ of global
descents of $u$, there is the chain $0<p_1<\dotsb< p_{k-1}<n$. In any other
chain there must be at least one equality, say $q_i=q_{i+1}$.
Then $\st(u_{q_i+1},\dotsc,u_{q_{i+1}})$ is the empty permutation and the
corresponding term is just the identity $1$, which is annihilated by
$\varphi$. 
Therefore, by~\eqref{E:cofree}, $\hat{\varphi}$ is given by
 \[\hat{\varphi}(\calM_u)\ =\
         \calM_{\st(u_1,\dotsc,u_{p_1})}\ten\dotsm\ten
      \calM_{\st(u_{p_{k-1}+1},\dotsc,u_{n})}\ \in\ V^{\ten k}\,.
 \]
Consider the map $\psi:V^{\ten k}\to(\SSym)^k$ that sends
\[ \calM_{v_{(1)}}\ten\dotsm\ten\calM_{v_{(k)}}\mapsto
\calM_{\zeta_\setT\cdot(v_{(1)}\times\dotsm\times v_{(k)})}\,,\]
where each
$v_{(i)}\in\frakS_{q_i}$ and
$\setT=\{q_1,q_1+q_2,\ldots,q_1+\cdots+q_{k-1}\}\subseteq[n{-}1]$.

Lemma~\ref{L:over/under} implies that
 $\GDes\bigl(\zeta_\setT\cdot(v_{(1)}\times\dotsm\times v_{(k)})\bigr)=\setT$,
since each $v_{(i)}$ has no global descents. Together with~\eqref{E:u^S} this
shows that $\hat{\varphi}\circ\psi=\id$.

On the other hand, letting $\setS=\GDes(u)$, Lemma~\ref{L:higherglobal} implies
that
 \[
  u=\zeta_\setS\cdot\bigl( \st(u_1,\dotsc,u_{p_1})\times\dotsm\times
      \st(u_{p_{k-1}+1},\dotsc,u_{n})\bigr)\,.
 \]
This shows that $\psi\circ\hat{\varphi}=\id$. Thus $\hat{\varphi}$ is an
isomorphism of graded coalgebras.
\end{proof}

\begin{rem} \label{R:poirier} If $V$ is finite dimensional then the graded dual
of $Q(V)$ is simply the (free) tensor algebra $T(V^*)$. More generally,
suppose $V=\oplus_{n\geq 1}V_n$ is a graded vector space for which
each component $V_n$ is finite dimensional. Then $Q(V)$ admits another grading,
for which the elements of $V_{n_1}\ten\dotsb\ten V_{n_k}$ have degree
$n_1+\cdots+n_k$ (with respect to the other grading, these elements
 have degree $k$).
 With respect to this new grading, the homogeneous components
are finite dimensional, and the graded dual of $Q(V)$ is
the tensor algebra on the graded dual of $V$ (again a free algebra).

In our situation, $\SSym=Q(V)$, with $V$ graded by the size
$n$ of the indexing permutations $u\in\frakS_n$.  
The corresponding grading on $\SSym$ is the original one,
for which $\calM_u$ has degree $n$ if $u\in\frakS_n$.
Its graded dual is therefore a free algebra. It is known that $\SSym$ is
self-dual with respect to this grading (see Section~\ref{S:duality}).
It follows that $\SSym$ is also a free algebra. This is a result of Poirier and
Reutenauer~\cite{PR95} who construct a different set of
algebra generators, not directly related to the monomial basis.
(See Remark~\ref{R:comtet}.)
\end{rem}

Let $C$ be a graded connected coalgebra.
The coradical $C^{(0)}$ of $C$ is the 1-dimensional component in degree 0
(identified with the base field via the counit).
The primitive elements of $C$ are
 \[
   \text{P}(C)\ :=\ \{x\in C\mid \Delta(x)=x\ten 1+1\ten x\}\,.
 \]
Set $C^{(1)}:=C^{(0)}\oplus \text{P}(C)$, the first level of the
coradical filtration.
More generally, the $k$-th level of the coradical filtration is
 \[
   C^{(k)}\ :=\ \bigl(\Delta^{(k)}\bigr)^{-1}
        \Bigl(\sum_{i+j=k}C^{\ten i}\ten C^{(0)}\ten C^{\ten j}\Bigl)\,.
 \]
We have $C^{(0)}\subseteq C^{(1)}\subseteq C^{(2)}
          \subseteq\dotsb\subseteq C=\bigcup_{k\geq 0}C^{(k)}$,
and
 \[
   \Delta(C^{(k)})\ \subseteq\ \sum_{i+j=k}C^{(i)}\ten C^{(j)}\,.
 \]
Thus, the coradical filtration measures the complexity of
iterated coproducts.

Suppose now that $C$ is a cofree graded coalgebra $Q(V)$.
Then the space of primitive elements is just $V$, and the $k$-th level
of the coradical filtration is $\oplus_{i=0}^k V^{\ten i}$. These are
straightforward consequences of the definition of the deconcatenation
coproduct.

Define
\[ \frakS_n^{(k)}\ :=\ \coprod_{i=0}^k \frakS_n^k \text{ \ and \ }
 \frakS^{(k)}\ :=\ \coprod_{i=0}^k \frakS^k \,.\]
 In other words,  $\frakS^{(0)}=\frakS_0$ and for
$k\geq 1$,
\[\frakS_n^{(k)}=\{u\in\frakS_n\mid u
                \text{ has at most $k{-}1$ global descents}\}\,.\]

 In Proposition~\ref{P:galoisglobal} we showed that
 $\GDes\colon\frakS_n\to\calQ_n$ is order-preserving.
 Since $\calQ_n$ is ranked by the cardinality of a subset, it follows that
  $\frakS_n^{(k)}$ is a lower order ideal of $\frakS_n$, with
  $\frakS_n^{(k)}\subseteq\frakS_n^{(k+1)}$.
 The coradical filtration corresponds precisely to this filtration of 
 the weak order on the symmetric groups by lower ideals.

 \begin{coro}\label{C:coradical}
 A linear basis for the $k$-th level of the coradical filtration of $\SSym$ is
 \[
   \{\calM_u\mid u\in \frakS^{(k)}\}\,.
 \]
 In particular, a linear basis for the space of primitive elements is
 \[
    \{\calM_u\mid u \text{ has no global descents}\}\,.
 \]
\end{coro}
\begin{proof} This follows from the preceding discussion.
\end{proof}

The original grading of $\SSym=\oplus_n\Q\frakS_n$ yields a grading on the
subspace $P(\SSym)$ of primitive elements and on each $(\SSym)^k$. 
Let $G_1(t)$ denote the Hilbert series of the space of primitive elements, or
equivalently, the generating function for the
set of permutations in $\frakS_n$ with no global descents,
 \[
   G_1(t)\ :=\ \sum_{n\geq 1}
           \dim_\Q\bigl(\text{P}_n(\SSym)\bigr)\, t^n\,.
 \]

More generally, let $G_k(t)$ be the Hilbert series of $(\SSym)^k$, or
equivalently, the generating function for permutations in $\frakS_n$ with
exactly $k-1$ global descents, 
 \[
   G_k(t)\ :=\ \sum_{n\geq k}
           \dim_\Q \bigl((\SSym)_n^k\bigr)\, t^n\,.
 \]
For instance,
 \begin{align*}
   G_1(t) & =t+t^2+  3t^3 +13t^4+71t^5+461t^6+3447t^7+\dotsb\\
   G_2(t) & =t^2+  2t^3 +7t^4+32t^5+177t^6+1142t^7+\dotsb\\
   G_3(t) & =t^3+ 3t^4 +12t^5+58t^6+327t^7+2109t^8+\dotsb\\
\end{align*}

There are well-known relationships  between the Hilbert
series of a graded space $V$, its powers $V^{\ten k}$ and their sum $Q(V)$.
In our case, these give the following formulas.

\begin{coro}\label{C:dimprimitive}
 We have
 \begin{enumerate}

  \item[({\it i})]
             ${\displaystyle
               \dim_\Q\bigl(\text{\rm P}_n(\SSym)\bigr)\ =\
        (-1)^{n-1}\begin{vmatrix}
               1! & 2! & \ldots &\ldots & n!\\
               1 & 1! & \ldots &\ldots & (n{-}1)!\\
               0 & 1 & 1! & \ldots & (n{-}2)!\\
               \vdots & \ddots & \ddots &\ddots&\vdots\\
               0 &\ldots& 0 & 1 & 1!
           \end{vmatrix}\ }$.

  \item[({\it ii})]${\displaystyle G_1(t)\ =\
                1-\frac{1}{\sum_{n\geq 0}n!\,t^n}}$\,.

  \item[({\it iii})]
          ${\displaystyle G_k(t)\ =\ \bigl(G_1(t)\bigr)^k}$\,.\rule{0pt}{18pt}
 \end{enumerate}
\end{coro}

\begin{rem}\label{R:comtet} Formula $(i)$ is analogous to a
 formula for ordinary descents in~\cite[Example~2.2.4]{St86}.
Formulas $(ii)$ and $(iii)$ in Corollary~\ref{C:dimprimitive} are 
due to Lentin~\cite[Section 6.3]{Len72}, see also
Comtet~\cite[Exercise VI.14]{Co74}. These references
do not consider global descents, but rather the problem of
decomposing a permutation $u\in\frakS_n$ as a non-trivial product
$u=v\times w$.
  This is equivalent to our study of global descents, as we may write
 $u=v\times w$ with $v\in\frakS_p$ exactly when $n{+}1{-}p$ is a global
 descent of $u\omega_n$.
 For instance, $u=56\mspace{2mu}324\mspace{1.mu}1$ has global descents
 $\{2,5\}$ and
 $u\omega_6=1\mspace{1.mu}423\mspace{2mu}65=1\times 312\times 21$.
  See the Encyclopedia of Integer Sequences~\cite{SlSeek} (A003319 and 
  A059438) for additional references in this connection.

 Poirier and Reutenauer~\cite{PR95} showed that the elements of the dual basis
 $\{\calF^*_u\}$ indexed by the connected permutations freely generate
 $(\SSym)^*$.
 Duchamp, Hivert, and Thibon dualize the resulting linear basis, giving a different
 basis than we do 
 for the space of  primitive 
 elements~\cite[Prop.~3.6]{DHT01}.
\end{rem}

\section{The descent map to quasi-symmetric functions}\label{S:descentmap}

We study the effect of the morphism of Hopf algebras~\eqref{E:descentmap}
 \[
   \calD\ \colon\ \SSym\ \to\ \QSym, \quad\text{ defined by }\quad
   \calF_u\ \mapsto\  F_{\Des(u)}\,
 \]
on the monomial basis.
Here, we use subsets $\setS$ of $[n{-}1]$ to index monomial and fundamental
quasi-symmetric functions of degree $n$, as discussed at the end of
Section~\ref{S:hopfquasi}.
Our main tool is the Galois connection $\frakS_n\rightleftarrows\calQ_n$ of
Section~\ref{S:descents}.

When we have a Galois connection between posets $P$ and $Q$ given by a
pair of  maps $f:P\to Q$ and $g:Q\to P$ as in~\eqref{E:galois},
a classical theorem of Rota~\cite[Theorem 1]{Rot} states that the M\"obius
functions of $P$ and $Q$ are related by
 \begin{equation*}
   \forall\ x\in P \text{ and }\ w\in Q,\quad \ 
      \sumsub{ y\in P\\x\leq y,\, f(y)=w}\!\!\mu_P(x,y)\ =\ 
      \sumsub{v\in Q\\v\leq w,\, g(v)=x}\!\!\mu_Q(v,w)\, .
 \end{equation*}
A conceptual proof of this simple but extremely useful result can be found
in~\cite{AF00}.

\begin{defi}\label{D:closed}
A permutation $u\in\frakS_n$ is {\em closed} if it is of the form
$u=\zeta_\setT$ for some $\setT\in\calQ_n$.
\end{defi}
Equivalently, in view of~\eqref{E:galoisdes} and~\eqref{E:galoisgdes}, $u$ is
closed if and only if $\Des(u)=\GDes(u)$.

{}From Proposition \ref{P:galois}, we deduce the
following fact about the M\"obius function of the weak order.

\begin{coro} 
 Let $u\in\frakS_n$ and $\setS\in\calQ_n$. Then
 \begin{equation}\label{E:mobius}
   \sumsub{u\leq v\in\frakS_n\\ \Des(v)=\setS}\mu_{\frakS_n}(u,v)\ =\ 
    \begin{cases}
         \mu_{\calQ_n}(\Des(u),\setS) & \text{if $u$ is closed,}\\
          0 & \text{if not.}
   \end{cases}
 \end{equation}
 \end{coro}

\begin{proof}
 Rota's formula says in this case that
 \[
    \sumsub{u\leq v\in\frakS_n\\\Des(v)=\setS\rule{0pt}{10pt}}
       \mu_{\frakS_n}(u,v)\ =\ 
    \sumsub{\setT\subseteq \setS\in\calQ_n\\\zeta_T=u\rule{0pt}{10pt}}
               \mu_{\calQ_n}(\setT,\setS)\,.
 \]
If $u$ is not closed, then the index set on the right hand side is empty. 
If $u$ is closed, then the index set consists only of the set $\setT=\Des(u)$,
by assertion (c) in the proof of Proposition \ref{P:galois}. 
\end{proof}

While there are explicit formulas for the M\"obius function of the weak order,
it is precisely the above result that allows us to obtain the description of
the map $\calD\colon \SSym\to\QSym$ in terms of the monomial bases.

\begin{thm}\label{T:map-monomial}
 Let $u\in\frakS_n$. Then
 \[
   \calD(\calM_u)\ =\ \begin{cases}
           M_{\GDes(u)} & \text{if $u$ is closed,}\\
                     0 & \text{if not.}
        \end{cases}
 \]
\end{thm}

\begin{proof}
 By definition, $\calM_u=\sum_{u\leq v} \mu_{\frakS_n}(u,v) \calF_v$, hence
 \begin{align*}
   \calD(\calM_u)
     &\ =\ \sum_{u\leq v} \mu_{\frakS_n}(u,v) F_{\Des(v)}\\
     &\ =\ \sum_{\setS}\Bigl(\sumsub{u\leq v\\\Des(v)=\setS}
                     \mu_{\frakS_n}(u,v) \Bigr)F_{\setS}\\
     &\ =\ \begin{cases}
   \sum_{\setS}\mu_{\calQ_n}(\Des(u),\setS)\,F_{\setS} & \text{if $u$ is closed}\\
          0 & \text{if not.}
\end{cases}
\end{align*}
 We complete the proof by noting that
 \[ 
   M_{\Des(u)}\ =\
   \sum_{\setS}\mu_{\calQ_n}(\Des(u),\setS)\,F_{\setS}\,
 \]
 by the definition of $M_{\Des(u)}$, and that since $u$ is closed, 
 $\Des(u)=\GDes(u)$.
\end{proof}

Malvenuto shows that  $\calD$  is a morphism of Hopf algebras by comparing the
structures on the fundamental bases of $\SSym$ and $\QSym$. 
We do the same for the monomial bases of $\SSym$ and
$\QSym$.

To compare the coproducts, first note that for any subsets
$\setS\subseteq[p-1]$ and $\setT\subseteq[q-1]$,
 \[
 \zeta_{\setS\cup\{p\}\cup\setT}=
 \zeta_{p,q}\cdot(\zeta_\setS\times\zeta_\setT)\,.
 \]
Therefore, if $u\in\frakS_n$ and $p\in\AGDes(u)$, then
 \[
    \text{$u$ is closed $\iff$ both $\st(u_1,\dotsc,u_p)$ and
                             $\st(u_{p+1},\dotsc,u_n)$ are closed.}
 \]
It follows that applying the map $\calD\colon \SSym\to\QSym$ to
formula~\eqref{E:cop-monomial} gives the usual
formula~\eqref{E:copqsym} for the coproduct of monomial quasi-symmetric
functions. 

For instance, we compare formula~\eqref{E:prod-M} with~\eqref{E:prod-ex}.
Since $\calD(\calM_{21})=M_{(1,1)}$ and $\calD(\calM_{12})=M_{(2)}$, applying
$\calD$ to~\eqref{E:prod-M} results in~\eqref{E:prod-ex}.
Indeed, the indices $u$ in the first row
of~\eqref{E:prod-M} all are closed, while none in the second row are closed.
It is easy to verify that the five terms on the right in the first row
in~\eqref{E:prod-M} map to the five terms on the right in~\eqref{E:prod-ex}.

The situation is different for the products. 
The geometric description of the
structure constants of the product on the monomial basis of
$\SSym$~\eqref{E:defalpha} admits an analogue for $\QSym$, but
this turns out to be very different from the known description in terms of
quasi-shuffles~\eqref{E:prodqsym}. 
We present this new description of the structure constants
for the product of monomial quasi-symmetric functions.

The role of the permutahedron is now played by the cube. 
Associating a subset $\setS$ of $[n{-}1]$ to its characteristic function
gives a bijection between subsets of $[n{-}1]$ and vertices of the 
$(n{-}1)$-dimensional cube $[0,1]^{n-1}$.
Coordinatewise comparison corresponds to subset inclusion, and the 1-skeleton
of the cube becomes the Hasse diagram of the Boolean poset $\calQ_n$.
In this way, we identify $\calQ_n$ with the vertices of the
$(n{-}1)$-dimensional cube. 

For each Grassmannian permutation $\zeta\in \Sh{(p,q)}$, consider the map
 \[
   r_\zeta\ \colon\ \calQ_p\times\calQ_q\ \to\ \calQ_{p+q}\,,
   \qquad (\setS,\setT)\ \mapsto\ 
   \Des((\zeta_\setS\times\zeta_\setT)\cdot\zeta^{-1})\,.
 \]
We describe this map $r_\zeta$ in more detail.
To that end, set 
 \[
   \mathrm{Cons}_p(\zeta)\ :=\ \{i\in[p{+}q{-}1]\mid 
                 \zeta^{-1}(i)+1=\zeta^{-1}(i{+}1)
                  \text{ and }  \zeta^{-1}(i)\neq p\}\,,
 \]
and recall that the vertices in a face of the cube are an interval in 
the Boolean poset, with every interval corresponding to a unique face.

\begin{lemm}\label{L:rzeta}
 Let $p,q$ be positive integers and $\zeta\in\Sh{(p,q)}$.
 The image of $r_\zeta$ is the face
 \[
    \bigl[\; \Des(\zeta^{-1}),\ 
      \Des(\zeta^{-1})\ \Disjoint\ \mathrm{Cons}_p(\zeta)
     \;\bigr]\,,
 \]
 which is isomorphic to the Boolean poset of subsets of\/ 
 $\mathrm{Cons}_p(\zeta)$.
\end{lemm}

\begin{proof}
 This is an immediate consequence of an alternative (and direct) description
 of $r_\zeta(\setS,\setT)$.
 For $\setT\in\calQ_q$, set $p+\setT:=\{p+t\mid t\in\setT\}$.
 Then, for $(\setS,\setT)\in\calQ_p\times\calQ_q$, we have
 \begin{equation}\label{E:rzeta}
    r_\zeta(\setS,\setT)\ =\ \Des(\zeta^{-1}){\textstyle \ \coprod\ }
       \Bigl(\mathrm{Cons}_p(\zeta)\cap
             \zeta\bigr(\setS\cup(p+\setT)\bigr)\Bigr)\,.
 \end{equation}
Assuming this for a moment, we note that
 the association $(\setS,\setT)\mapsto\zeta\bigl(\setS\cup(p+\setT)\bigr)$ 
 is a bijection between $\calQ_p\times\calQ_q$ and
 subsets of $\{i\mid \zeta^{-1}(i)\neq p\}$. Intersecting with
$\mathrm{Cons}_p(\zeta)$ we obtain a surjection onto subsets of
$\mathrm{Cons}_p(\zeta)$, which yields the desired
description of the image of $r_\zeta$.

 We prove~\eqref{E:rzeta}.
 Let $(\setS,\setT)\in\calQ_p\times\calQ_q$ and set 
 $w:=(\zeta_\setS\times\zeta_\setT)\cdot\zeta^{-1}$ so that
 $\Des(w)=r_\zeta(\setS,\setT)$.
 Note that $\Des(\zeta_\setS\times\zeta_\setT)=\setS\cup(p+\setT)$ (this is a
particular case of Lemma~\ref{L:over/under}) and  if $i\leq p<j$, then 
 $(\zeta_\setS\times\zeta_\setT)(i)\leq p<
  (\zeta_\setS\times\zeta_\setT)(j)$.

 Let $i\in[n{-}1]$.
 We consider whether or not $i$ is a descent of $w$.
 First, suppose $i\in\Des(\zeta^{-1})$.
 Since the values $1,2,\dotsc,p$ and $p{+}1,p{+}2,\dotsc,p{+}q$
occur in  order in the permutation $\zeta^{-1}$
(because $\zeta\in\Sh{(p,q)}$),  we must have
$\zeta^{-1}(i)>p\geq\zeta^{-1}(i{+}1)$ and so  $w(i)>p\geq w(i+1)$, thus
$i\in\Des(w)$.

 Now suppose that $i$ is not a descent of $\zeta^{-1}$.
 If $\zeta^{-1}(i)+1<\zeta^{-1}(i{+}1)$, then we must have 
 $\zeta^{-1}(i)\leq p < \zeta^{-1}(i{+}1)$, again because
 $\zeta\in\Sh{(p,q)}$. Hence $w(i)\leq p < w(i+1)$
and $i$ is not a descent of  $w$.
 If instead we have $\zeta^{-1}(i)+1=\zeta^{-1}(i{+}1)$, then there are two
 cases to consider.
 If $i=\zeta^{-1}(p)$, then this forces $\zeta$ to be $1_{p+q}$ so
 $w(i)=w(p)\leq p<w(i{+}1)$,  and we conclude that $i$ is not a descent of $w$.
 If $i\neq\zeta^{-1}(p)$, then $i\in\mathrm{Cons}_p(\zeta)$ and we see
 that $i$ is a descent of $w$ exactly when
 $\zeta^{-1}(i)\in\setS\cup(p+\setT)$.
 This proves~\eqref{E:rzeta} and completes the proof of the lemma.
\end{proof}

\begin{figure}[htb]
 \[
   \begin{picture}(175,135)(-11,8)
    \put(20,11){\epsfysize=120pt\epsfbox{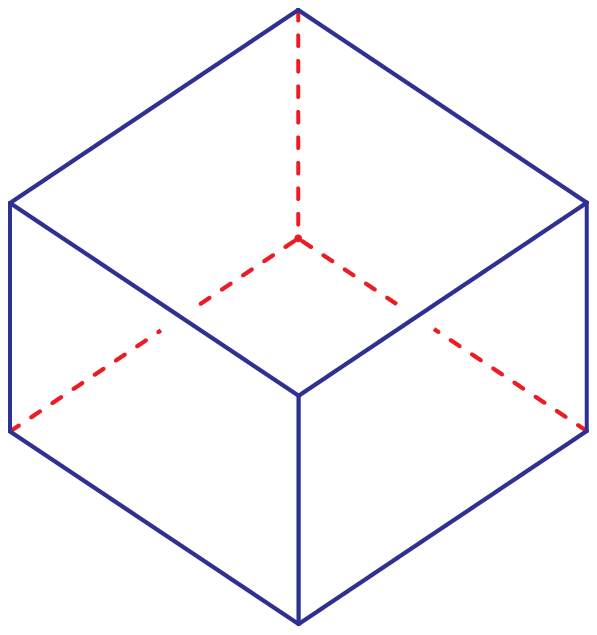}}
                            \put(57,138){$\{1,2,3\}$}
    \put(-11,91){$\{2,3\}$} \put(80, 90){$\{1,3\}$}  \put(136,91){$\{1,2\}$}
    \put(  0,46){$\{3\}$}   \put(55, 49){$\{2\}$}    \put(136,46){$\{1\}$}
                            \put(74, 0){$\emptyset$}
   \end{picture}
 \]
 \caption{Vertices of the cube}
\label{F:labelcube}
\end{figure}

\begin{figure}[htb]
 \[
  \begin{picture}(155,157)(10,-22)
   \put(40,15){\epsfysize=100pt\epsfbox{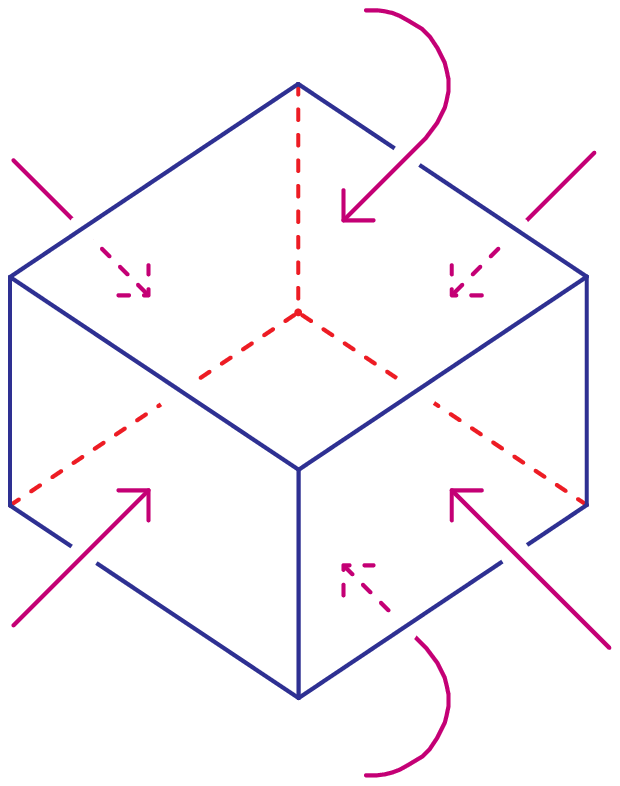}}
   \put(10,110){$r_{2341}(3)$}  \put(130,110){$r_{4123}(1)$}
   \put(57,125){$r_{3412}(2)$} 

   \put(10,15){$r_{1234}(1)$}   \put(130,15){$r_{1234}(3)$}
   \put(57, 0){$r_{1234}(2)$} 

   \put(80,-22){(a)}
  \end{picture}
  \qquad\qquad
  \begin{picture}(195,157)(-38,-25)
  \put( 0,15){\epsfysize=90pt\epsfbox{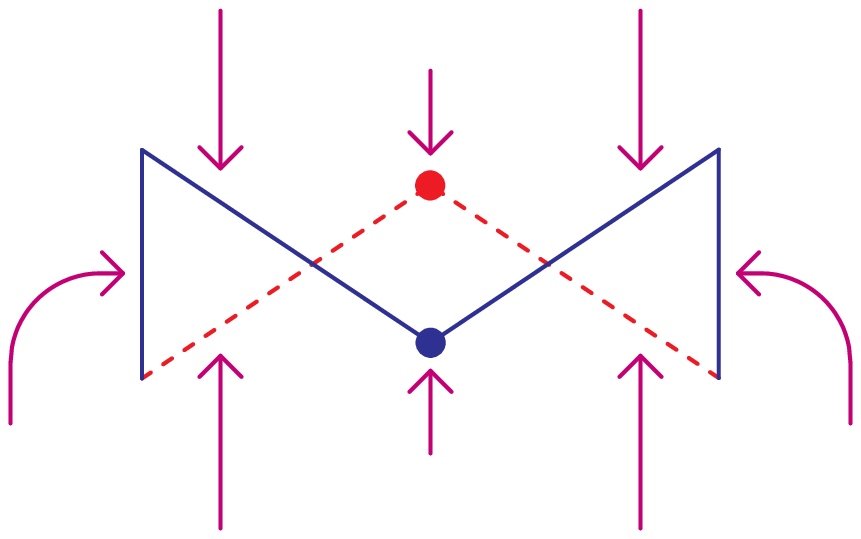}}
    \put(  2,114){$r_{1342}(3)$}  \put(82,114){$r_{3124}(1)$} 
           \put(41,103){$r_{2413}(2)$} 
    \put(  2,3){$r_{1243}(3)$}  \put(82,3){$r_{2134}(1)$} 
           \put(41,17){$r_{1324}(2)$} 

    \put(-38,22){$r_{1423}(2)$}  \put(121,22){$r_{2314}(2)$} 
   \put(48,-25){(b)}
   \end{picture}
 \]
\caption{(a) The  facets of the cube: $r_{1234}(\calQ_p\times\calQ_q)$ and
\(r_{\zeta_{p ,q}}(\calQ_p\times\calQ_q)\ (=r_{\zeta_{p ,q}}(p))\).
\mbox{\hspace{5.3em}}(b)  The edges and vertices
$r_\zeta(\calQ_p\times\calQ_q)$, $\zeta\neq 1234,\zeta_{p,q}.$}
\label{F:typescube} \end{figure}

Unlike the case of the permutahedron, the image of $r_\zeta$ need not be a
facet.
Indeed, by Lemma~\ref{L:rzeta}, the image of $r_\zeta$ is a
facet only if $\#\mathrm{Cons}_p(\zeta)=p+q-2$, and this occurs only when 
$\zeta=1_{p+q}$ or $\zeta=\zeta_{p,q}$.
Figure~\ref{F:labelcube} displays the vertices of the 3-cube and 
Figure~\ref{F:typescube} shows which faces occur as the image 
$r_\zeta(\calQ_p\times\calQ_q)$. 
Observe that while not all faces occur as images of some
$r_\zeta(\calQ_p\times\calQ_q)$, any face that does occur is the image of a
unique such map.
This is the general case.

\begin{lemm}
 A face of \/$\calQ_n$ is the image  of $\calQ_p\times\calQ_{n-p}$ under a map
 $r_\zeta$ for at most one pair $(\zeta,p)$.
\end{lemm}

\begin{proof}
 Suppose $\zeta\in\Sh{(p,n-p)}$ for some $0<p<n$.
 We will observe that the pair of sets $\Des(\zeta^{-1})$ and
 $\mathrm{Cons}_p(\zeta)$ determines $\zeta$ and $p$ uniquely by describing 
 these sets.

 Suppose first that $\zeta=1_n$.
 Then $\Des(\zeta^{-1})=\emptyset$ and 
 $\mathrm{Cons}_p(\zeta)=[n{-}1]-\{p\}$.

 Suppose now that $\zeta\in\Sh{(p,n-p)}$ is not the identity permutation.
 Then $\zeta$ determines $p$ and $\Des(\zeta^{-1})\neq\emptyset$.
 Since the values $1,2,\dotsc,p$ and
 $p{+}1,\dotsc,n$ occur in order in $\zeta^{-1}$,
 there exist numbers
 \[
    0\leq b_0 < a_1 < b_1< \dotsb < a_k < b_k \leq n
 \]
 such that the values in $[p]$ occur in order in the intervals
 \[
    [0,b_0], [a_1+1,b_1], \dotsc, [a_k+1, b_k]\,,
 \]
 and the values in $\{p{+}1,\dotsc,n\}$ in the complementary set.
 Thus $\Des(\zeta^{-1})=\{a_1,\ldots,a_k\}$ and 
 $\mathrm{Cons}_p(\zeta)=[n{-}1] - \{b_0,a_1,b_1,a_2,\ldots,a_k,b_k\}$.

 It follows that $\zeta$ and $p$ determine and are determined by the sets
 $\Des(\zeta^{-1})$ and $\mathrm{Cons}_p(\zeta)$, which completes the proof of
 the lemma.
\end{proof}

\begin{thm}\label{T:prodmonqsym} 
  Suppose $p,q$ are positive integers.
 Let $\setS\subseteq[p{-}1]$, $\setT\subseteq[q{-}1]$ and 
 $\setR\subseteq[p{+}q{-}1]$. 
 The coefficient of $M_{p+q,\setR}$ in $M_{p,\setS}\cdot M_{q,\setT}$ is
 \begin{equation}\label{E:defalphaqsym}
   \#\{\zeta\in \Sh{(p,q)}\mid (\setS,\setT)\ =\ 
    \max r_\zeta^{-1}[\emptyset,\setR]\}\,.
 \end{equation}
 In other words, this coefficient counts the
 number of faces of the cube of type $(p,q)$ with the property that the vertex
 $r_\zeta(\setS,\setT)$ is below $\setR$  and it is
 the maximal vertex in the the face $r_\zeta(\calQ_p\times\calQ_q)$ 
 below $\setR$.
\end{thm}

\begin{proof} 
 By Theorem~\ref{T:map-monomial}, 
 $M_{p,\setS}\cdot M_{q,\setT}=
   \calD(\calM_{\zeta_\setS}\cdot\calM_{\zeta_\setT})$.
 We expand the product using Theorem~\ref{T:prod-monomial}, and then apply
 the map $\calD$ and Theorem~\ref{T:map-monomial} to obtain
\[   M_{p,\setS}\cdot M_{q,\setT} 
    \ =\ \calD(\calM_{\zeta_\setS}\cdot\calM_{\zeta_\setT})
    \ =\ \calD\Bigl(\sum_{w\in\frakS_{p+q}} 
               \alpha^w_{\zeta_\setS, \zeta_\setT} \calM_w \Bigr)
    \ =\ \sum_{\setR\in\calQ_{p+q}}
             \alpha^{\zeta_R}_{\zeta_\setS,\zeta_\setT} M_{p+q,\setR}\,.
\]
 According to~\eqref{E:defalpha},
 \[
   \alpha^{\zeta_\setR}_{\zeta_\setS,\zeta_\setT}\ =\ 
    \#\{\zeta \in \Sh{(p,q)}\mid (\zeta_\setS,\zeta_\setT)\ =\ 
         \max\rho_\zeta^{-1}[1,\zeta_\setR]\}\,. 
 \]
 By Proposition~\ref{P:galois}, for any $\setS$, $\setT$, and $\setR$ we
 have
 \[
   \Des\bigl((\zeta_\setS\times\zeta_\setT)\cdot\zeta^{-1}\bigr)
   \ \subseteq\ \setR \quad\iff\quad
   (\zeta_\setS\times\zeta_\setT)\cdot\zeta^{-1}\ \leq\ \zeta_\setR\,.
 \]
 In other words, 
 \[
   r_\zeta(\setS,\setT)\ \leq\ \setR \quad\iff\quad
    \rho_\zeta(\zeta_\setS,\zeta_\setT)\ \leq\ \zeta_{\setR}\,.
 \]
 This implies that the structure constant
 $\alpha^{\zeta_R}_{\zeta_\setS,\zeta_\setT}$ is as stated.
\end{proof}

\begin{figure}[!htb]
 \[
  \begin{picture}(270,115)
   \put(38,1){\epsfysize=110pt\epsfbox{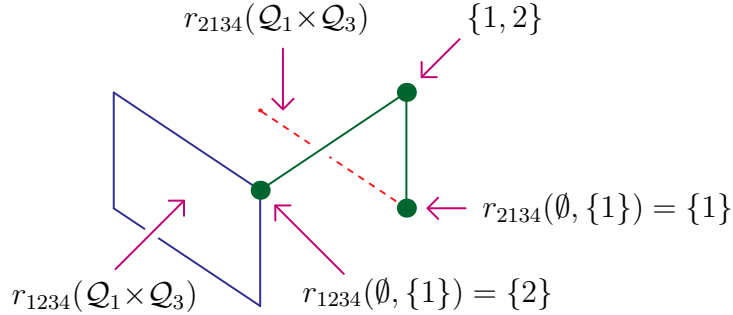}}
   \put(0,2){$r_{1234}(\calQ_1{\times}\calQ_3)$}
   \put(110,3){$r_{1234}(\emptyset,\{1\})=\{2\}$}
   \put(65,107){$r_{2134}(\calQ_1{\times}\calQ_3)$}
   \put(172,107){$\{1,2\}$}
   \put(178,35){$r_{2134}(\emptyset,\{1\})=\{1\}$}
  \end{picture}
 \]
  \caption{The faces  $r_{1234}$ and $r_{2134}$ of type $(1,3)$, and the vertex
$\{1,2\}$\label{F:alphacube}} \end{figure}

We give an example.
Let $p=1$, $q=3$, $\setS=\emptyset$ and $\setT=\{1\}$.  
In terms of compositions, we have $M_{\emptyset,1}=M_{(1)}$, and
     $M_{\{1\},3}=M_{(1,2)}$.
Equation~\eqref{E:prodqsym} gives
 \begin{align*}
  M_{\emptyset,1}\cdot M_{\{1\},3}
    \ =\ M_{(1)}\cdot M_{(1,2)}&\ =\ 
    2M_{(1,1,2)}+M_{(1,2,1)}+M_{(2,2)}+M_{(1,3)}\\
    &\ =\ 2M_{\{1,2\},4}+M_{\{1,3\},4}+M_{\{2\},4}+M_{\{1\},4}\,.
 \end{align*}
On the other hand, \eqref{E:defalphaqsym} also predicts that the coefficient of
$M_{\{1,2\}}$ is $2$. 
Of the four possible faces of type $(1,3)$, only two satisfy the required
condition.  
One corresponds to the shuffle $1234$ (it is a facet) and the other to $2134$
(it is an edge). 
They are shown in Figure~\ref{F:alphacube}, together with the vertices
$r_{1234}(\emptyset,\{1\})=\{2\}$, 
$r_{2134}(\emptyset,\{1\})=\{1\}$, and the vertex $\{1,2\}$.

\section{$\SSym$ is a crossed product over $\QSym$}\label{S:crossed}

We obtain a decomposition of the algebra
structure of $\SSym$ as a crossed product over the Hopf algebra $\QSym$.
We refer the reader to~\cite[\S 7]{Mo93a} for a review of this construction
in the general Hopf algebraic setting.  Let us only say that the crossed
product of a Hopf algebra $K$ with an algebra $A$ with respect to a Hopf cocycle
$\sigma:K\ten K\to A$ is a certain algebra structure on the space $A\otimes K$,
denoted by $A\#_{\sigma} K$.

\begin{thm}
  The map $\calZ\colon\QSym\to\SSym$, $M_{\setS}\mapsto \calM_{\zeta_\setS}$,
 is a morphism of  coalgebras and a right inverse
 to the morphism of Hopf algebras $\calD\colon\SSym\to\QSym$.
\end{thm}

\begin{proof} This is
immediate from Theorems~\ref{T:cop-monomial} and~\ref{T:map-monomial}.
\end{proof}

In this situation, an important theorem of Blattner, Cohen, and
Montgomery~~\cite{BCM} applies. Namely, suppose $\pi:H\to K$ is a morphism of
Hopf algebras that admits a coalgebra splitting (right inverse) $\gamma:K\to H$. Then
there is a \emph{crossed product} decomposition
\[H\ \cong\ A\#_{\sigma}K\]
where $A$, a subalgebra of $H$, is the \emph{left Hopf kernel} of $\pi$:
\[A\ =\ \{h\in H \mid \sum h_1\ten\pi(h_2)=h\ten 1\}\]
and the \emph{Hopf cocycle} $\sigma:K\ten K\to A$ is
\begin{equation}\label{E:cocycle}
\sigma(k,k')=\sum\gamma(k_1)\gamma(k'_1)S\gamma(k_2k'_2)\,.
\end{equation}
This result, as well as some generalizations,
can be found in~\cite[\S 7]{Mo93a}. Note that if $\pi$ and $\gamma$ preserve
gradings, then so does the rest of the structure.

Let $A$ be the left Hopf kernel of $\calD\colon\SSym\to\QSym$ and
$A_n$ its $n$-th homogeneous component.
Once again the monomial basis of $\SSym$ proves useful in
describing $A$.

\begin{thm} A basis for  $A_n$ is the
set  $\{\calM_u\}$ where $u$ runs over all permutations of $n$ that
are not of the form
\begin{equation}\label{E:hopfkernel}\tag{$*$}
*\ldots*12\ldots n{-}k
\end{equation}
for any $k=0,\ldots,n-1$. In particular,
\[\dim A_n\ =\ n!-\sum_{k=0}^{n-1}k!\,.\]
\end{thm}

\begin{proof} By the theorem of Blattner, Cohen, and Montgomery,
 $\SSym\cong A\#_{\sigma}\QSym$, in particular
 $\SSym\cong A\ten\QSym$ as vector spaces.
 The generating functions for the dimensions of these
 algebras are therefore related by
 \[\sum_{n\geq 0}^\infty n!t^n\ =\ \sum_{n\geq 0} a_nt^n\cdot
 \Bigl(1+\sum_{n\geq 1}2^{n-1}t^n\Bigr)\ = \ 
\sum_{n\geq 0} a_nt^n\cdot \frac{1}{1-\sum_{n\geq 1}t^n} \,.\]
 It follows that $a_n= n!-\sum_{k=0}^{n-1}k!$ as claimed.

Observe that $a_n$ counts the 
permutations in $\frakS_n$ that are not of
the form \eqref{E:hopfkernel}. Since the $\calM_u$ are linearly
independent,  it suffices to show that if $u$ is not of that form
 then $\calM_u$ is  in the Hopf kernel. 
 Now, for any $u\in\frakS_n$ and $p\in\GDes(u)$, 
 we have that $\st(u_{p+1},\ldots,u_n)=(u_{p+1},\ldots,u_n)$.
Hence, if $u$ is not of the form  \eqref{E:hopfkernel}, the same is true of
$\st(u_{p+1},\ldots,u_n)$ and therefore this permutation is not closed.
It follows from
Theorems~\ref{T:cop-monomial} and~\ref{T:map-monomial} that
$(id\ten\calD)\Delta(\calM_u)=\calM_u\ten 1$.
\end{proof}

 \begin{rem}  These results were motivated by a question of Nantel Bergeron, who
 asked (in dual form) if $\SSym$ is cofree as right comodule over $\QSym$.
 This is an immediate consequence of the crossed product decomposition.
\end{rem}
 
Consider again the general situation of a morphism of Hopf algebras
$\pi:H\to K$ with a coalgebra splitting $\gamma:K\to H$. 
This induces an exact sequence of Lie algebras
\begin{equation}\label{E:lie-extension}
 0\to \textrm{P}(H)\cap A\to \textrm{P}(H)\map{\pi} \textrm{P}(K)\to 0
\end{equation}
with a linear splitting $\textrm{P}(K)\map{\gamma}\textrm{P}(H)$, 
where $\textrm{P}(H)$ denotes the space
of primitive elements of $H$, viewed as a Lie algebra under the commutator
bracket $\left[h,h'\right]=hh'-h'h$.

The Hopf cocycle restricts to a linear map
$\sigma:\textrm{P}(K)\ten \textrm{P}(K)\to \textrm{P}(H)\cap A$; 
in fact, for primitive elements $k$ and
$k'$, \eqref{E:cocycle} specializes to
\begin{equation}\label{E:cocycle-on-prim}
\sigma(k,k')=S\gamma(kk')-\gamma(k')\gamma(k)
\end{equation}
and a direct calculation shows that this element of $H$ is primitive. Moreover,
the Lie cocycle corresponding to \eqref{E:lie-extension} is the map
$\tilde{\sigma}:\textrm{P}(K)\wedge \textrm{P}(K)\to \textrm{P}(H)\cap A$ given by
\begin{equation}\label{E:lie-cocycle}
\tilde{\sigma}(k,k')=\left[\gamma(k),\gamma(k')\right]-\gamma(\left[k,k'\right])
= \sigma(k,k')-\sigma(k',k) \,.
\end{equation}
This map is a non-abelian Lie cocycle in the sense that the following conditions
hold. For $k$, $k'\in \textrm{P}(K)$ and $a\in \textrm{P}(H)\cap A$,
\begin{gather*}
k\cdot(k'\cdot a)-k'\cdot(k\cdot
a)=\left[\tilde{\sigma}(k,k'),a\right]+\left[k,k'\right]\cdot a \\
k\cdot\tilde{\sigma}(k',k'')-k'\cdot\tilde{\sigma}(k,k'')+
k''\cdot\tilde{\sigma}(k,k')=\tilde{\sigma}(\left[k,k'\right],k'')-
\tilde{\sigma}(\left[k,k''\right],k')+\tilde{\sigma}(\left[k',k''\right],k)
 \end{gather*}
 where $k\cdot a=\left[\gamma(k),a\right]$.

Let us apply these considerations to the morphism  $\SSym\map{\calD}\QSym$ and
the coalgebra splitting $\QSym\map{\calZ}\QSym$.
The structure constants of the Hopf cocycle $\sigma$ do not have constant sign.
However, its restriction to primitive elements of $\QSym$ has non-negative
structure constants on the monomial bases. They turn out to be
particular structure constants  of the product of $\SSym$. 

Recall that these structure constants $\alpha^w_{u,v}$ are defined for
$u\in\frakS_p$, $v\in\frakS_q$ and $w\in\frakS_{p+q}$ by the identity
 \[
    \calM_u\cdot\calM_v\ = \ \sum_{w\in\frakS_{p+q}}\alpha^w_{u,v}\calM_w\,.
 \]
The combinatorial description of these constants showing their non-negativity
is given by~\eqref{E:A-def}.

\begin{lemm}\label{L:struc-constants}
 For $p,q\geq 1$, and $w\in \frakS_{p+q}$ closed, we have
 $\alpha^w_{1_p,1_q}=0$ except in the following cases
 \[
   \alpha_{1_p,1_q}^{1_{p+q}}\ =\ 1\,,\ 
   \alpha_{1_p,1_p}^{\zeta_{p,p}}\ =\ 2\quad 
    \textrm{\ and if $p\neq q$, then\ }\quad 
   \alpha_{1_p,1_q}^{\zeta_{p,q}}\ =\ 1\,.
 \]
\end{lemm}

\begin{proof}
 Apply the map $\calD$ to the product
 \[
   \sum_{w\in\frakS_{p+q}}\alpha^w_{1_p,1_q}\calM_w\ =\ 
    \calM_{1_p}\cdot\calM_{1_q}\,,
 \]
 to obtain (using \eqref{E:prodqsym})
 \[
    \sum_{w\in\frakS_{p+q}}\alpha^w_{1_p,1_q}\calD(\calM_w)\ =\ 
    M_{(p)}\cdot M_{(q)}\ =\ M_{(p,q)} + M_{(q,p)} + M_{(p+q)}\,.
 \]
 The result is immediate, as $\calD(\calM_w)=0$ unless $w$ is closed,
 and  we have $\calD(\calM_{\zeta_{p,q}})=M_{(p,q)}$ and 
 $\calD(\calM_{1_{p+q}})=M_{(p+q)}$.
\end{proof}

We use this lemma to give a
combinatorial description of $\sigma$ and the Lie
cocycle $\tilde{\sigma}$ on primitive elements.
By~\eqref{E:copqsym}, $\{M_{(n)}\}_{n\geq 1}$ is a linear basis 
for the space of primitive elements of $\QSym$.
Thus $\textrm{P}(\QSym)$ is an abelian Lie algebra with each homogeneous
component of dimension $1$.
Recall that $\{\calM_u\mid\textrm{$u$ has no global descents}\}$ is a basis of
the primitive elements of $\SSym$, and thus $A\cap\textrm{P}(\SSym)$ has a basis given by
those $\calM_u$ where $u$ has no global descents and $u$ is not an identity
permutation, $1_n$.  

\begin{thm}\label{P:cocycle}
  For any $p$, $q\geq 1$,
 \begin{align*}
  \sigma(M_{(p)},M_{(q)})&=\sum_{w\neq
  \zeta_{p,q},\,\zeta_{q,p},\,1_{p+q}}\alpha_{1_q,1_p}^{w}\calM_w \\
  \tilde{\sigma}(M_{(p)},M_{(q)})&=\sum_{w}\bigl(\alpha_{1_q,1_p}^{w}-
  \alpha_{1_p,1_q}^{w}\bigr)\calM_w\,.
 \end{align*}
\end{thm}

\begin{proof} 
 Since $M_{(p)}\cdot M_{(q)}=M_{(p,q)}+M_{(q,p)}+M_{(p+q)}$,~\eqref{E:cocycle-on-prim} gives 
 \begin{align*}
  \sigma(M_{(p)},M_{(q)})&\ =\ S\calZ(M_{(p,q)}+M_{(q,p)}+M_{(p+q)})-
  \calZ(M_{(q)})\cdot \calZ(M_{(p)}\\
  &\ =\ 
 S(\calM_{\zeta_{p,q}}+\calM_{\zeta_{q,p}}+\calM_{1_{p+q}})-\calM_{1_q}\cdot
  \calM_{1_p}\,.
 \end{align*}
Using~\eqref{E:antipode} and~\eqref{E:cop-monomial},
we compute
$S(\calM_{\zeta_{p,q}})=\calM_{1_p}\cdot\calM_{1_q}-\calM_{\zeta_{p,q}}$
and $S(\calM_{1_{p+q}})=-\calM_{1_{p+q}}$.  
Therefore,
 \[
   \sigma(M_{(p)},M_{(q)})\ =\  \calM_{1_p}\cdot\calM_{1_q}-\calM_{\zeta_{p,q}}
    -\calM_{\zeta_{q,p}}-\calM_{1_{p+q}}\,.
 \]

The formula for $\sigma(M_{(p)},M_{(q)})$ follows by expanding the product
and using Lemma~\ref{L:struc-constants}.
The expression for $\tilde{\sigma}$ follows immediately
from~\eqref{E:lie-cocycle}. 
\end{proof}

\section{Self-duality of $\SSym$ and applications}\label{S:duality}

The Hopf algebra $\SSym$ is self-dual. This appears 
in~\cite[section 5.2]{Malv}, ~\cite[Theorem 3.3]{MR95},
and~\cite{Jol99}. We provide a proof below, for completeness. 
We investigate the combinatorial implications of 
this self-duality, particularly when expressed in terms of the monomial basis.
We explain how a result of
Foata and Sch\"utzenberger on the numbers
 \[
   d(\setS,\setT)\ =\ \#\{x\in\frakS_n\mid \Des(x)=\setS,\
   \Des(x^{-1})=\setT\}
 \]
is a consequence of self-duality of $\SSym$ and obtain analogous results
for the numbers
 \[
   \theta(u,v)\ :=\ \#\{x\in\frakS_n\mid x\leq u,\ x^{-1}\leq v\}\,. 
 \]

The Hopf algebra $\SSym$ is connected and graded with each homogeneous
component finite dimensional. We consider its \emph{graded dual} $(\SSym)^*$
whose homogeneous component in degree $n$ is the linear dual of  the
homogeneous component in degree $n$ of $\SSym$. Let $\{\calF^*_u\mid
u\in\frakS_n,\ n\geq 0\}$ and $\{\calM^*_u\mid u\in\frakS_n,\ n\geq 0\}$ be
the bases of $(\SSym)^*$ dual to the fundamental and monomial bases of $\SSym$,
respectively. $(\SSym)^*$ is another graded connected Hopf algebra.

\begin{thm} \label{T:selfdual} 
  The map
 \begin{equation}\label{E:selffund}
   \Theta\colon (\SSym)^*\to\SSym\,,\qquad \calF^*_u\mapsto\calF_{u^{-1}}
\end{equation}
is an isomorphism of Hopf algebras.
On the monomial basis it is given by
 \begin{equation}\label{E:selfmon}
   \Theta(\calM^*_u)\ =\ \sum_v\theta(u,v)\calM_{v}\,.
 \end{equation}
\end{thm}

\noindent{\it Proof.} 
 Note that $\Theta^*=\Theta$.
 Therefore, it suffices to show that $\Theta$ is a morphism of coalgebras.
 We rewrite the product~\eqref{E:prod-fundamental} of $\SSym$.
 Let $u\in\frakS_p$ and $v\in\frakS_q$.
 Then
 \[
   \calF_u\cdot\calF_v\ =\ \sum_{w\in\frakS_{p+q}}
        \#\{\zeta\in\Sh{(p,q)}\mid (u\times v)\cdot\zeta^{-1}=w\}\,\calF_w\,.
 \]
 Therefore the  (dual) coproduct of $(\SSym)^*$ is
 \[
   \Delta(\calF_w^*)\ =\ \sum_{p+q=n}\ 
           \sum_{u\in\frakS_p,v\in\frakS_q}
    \#\{\zeta\in\Sh{(p,q)}\mid (u\times v)\cdot\zeta^{-1}=w\}\;
       \calF_u^*\ten \calF_v^*\,.
 \]
 On the other hand, as observed in~\ref{R:loday}, the coproduct of $\SSym$  
 can be written as
 \[
   \Delta(\calF_w)\ =\ \sum_{p+q=n}\
         \sum_{u\in\frakS_p,v\in\frakS_q}
          \#\{\zeta\in\Sh{(p,q)}\mid \zeta\cdot(u\times v)\:=\:w\}\;
         \calF_u\ten \calF_v\,.
 \]
 It follows that $\Theta$ is a morphism of coalgebras because
 \[
   w\ =\ \zeta\cdot(u\times v)\ \ \iff\ \ 
     w^{-1}\ =\ (u^{-1}\times v^{-1})\cdot\zeta^{-1}\,.
 \]

 Since $\calF_u=\sum_{u\leq x}\calM_x$, we have
 $\calM_u^* =\sum_{x\leq u}\calF_x^*$.
 Therefore,
 \[
   \Theta(\calM_u^*)\ =\ \sum_{x\leq u}\calF_{x^{-1}}
      \ =\ \sum_{x\leq u}\:\sum_{x^{-1}\leq v}\calM_v
      \ =\ \sum_v \theta(u,v)\calM_v\,. 
    \eqno{\Box}\raisebox{-30pt}{\rule{0pt}{0pt}}
 \]

Formula~\eqref{E:selfmon} for the morphism $\Theta$ of Hopf algebras has
combinatorial implications which we develop.
Recall that $\alpha^w(u,v)$ and $\kappa(u,w)$ denote the structure constants
of the product and antipode of $\SSym$ in terms of the monomial basis.
These integers were described in
Theorems~\ref{T:prod-monomial} and~\ref{T:ant-monomial}. 
Consider $\theta$, $\alpha^w$, and $\kappa$ to be 
matrices with rows and columns indexed by elements of $\frakS_n$.

\begin{thm}\label{T:theta-form}
  For any $u\in\frakS_p$, $v\in\frakS_q$, and $w\in\frakS_{p+q}$, we have
 \begin{itemize}
  \item[({\it i})] \qquad 
   ${\displaystyle
       (\theta\alpha^w\theta)(u,v)\ =\ 
           \theta(\zeta_{p,q}{\cdot}(u\times v),\,w)}$\,, 
  \item[({\it ii})]\qquad 
       $\kappa^t\theta\ =\ \theta\kappa$.  
 \end{itemize}
\end{thm}

\begin{proof} 
 By Lemma~\ref{L:global}, the coproduct of $\SSym$~\eqref{E:cop-monomial} can
 be written as 
 \[
   \Delta(\calM_w)\ =\ \sum_{p+q=n} \ 
         \sumsub{u\in\frakS_p,v\in\frakS_q\\\zeta_{p,q}\cdot(u\times v)=w} 
    \calM_u\ten \calM_v\,,
 \]
 Therefore, the dual product is
 \[
    \calM_u^*\cdot\calM_v^*\ =\ \calM_{\zeta_{p,q}\cdot(u\times v)}^*\,.
 \]
 Thus, the right hand side of~({\it i}) is the coefficient of $\calM_w$ in
 $\Theta(\calM_u^*\cdot\calM_v^*)$. 
 On the other hand, since $\theta(u,v)=\theta(v,u)$, we have
 \[
   (\theta\alpha^w\theta)(u,v)\ =\ 
    \sum_{x,y\in\frakS_n}\theta(u,x)\alpha^w(x,y)\theta(y,v)\ =\ 
               \sum_{x,y\in\frakS_n}\theta(u,x)\theta(v,y)\alpha^w(x,y)\,.
 \]
 Thus the left hand side of~({\it i}) is the coefficient of $\calM_w$ in
 $\Theta(\calM_u^*)\cdot\Theta(\calM_v^*)$.  
 Since $\Theta$ is a morphism of
 algebras,~({\it i}) holds.

 The second formula directly expresses that $\Theta$ preserves
 antipodes, since the antipode of $(\SSym)^*$ is the dual of the antipode of
 $\SSym$. 
\end{proof}

One may view Theorem~\ref{T:theta-form}({\it i}) as a recursion reducing the
computation of $\theta(u,v)$ to the case when $u$ and $v$ have no global
descents, by virtue of Lemma~\ref{L:global}. 
On the other  hand, since $\theta$ is an invertible  matrix, this 
 provides another  semi-combinatorial description of the structure
constants $\alpha^w(u,v)$. 

One may also impose the condition that $\Theta$ preserves coproducts, but this
leads again to ({\it i}) of Theorem~\ref{T:theta-form}.
On the other hand, the equivalent of  ({\it ii}) of Theorem~\ref{T:theta-form}
for the fundamental basis leads to the following non-trivial identity.

\begin{prop} 
 For any $u$ and $v\in\frakS_n$,
 \begin{multline*}
   \#\{\setS\subseteq[n{-}1]\mid \Des(vu_\setS)\subseteq\setS 
               \text{ and $\#\setS$ is odd}\}\\
    \shoveleft{+ \#\{\setS\subseteq[n{-}1]\mid\Des(uv_\setS)\subseteq\setS 
                  \text{ and $\#\setS$ is even}\}} \\
   \shoveright{=\ \#\{\setS\subseteq[n{-}1]\mid \Des(vu_\setS)\subseteq\setS
   \text{ and $\#\setS$ is even}\}\quad}\\
    + \#\{\setS\subseteq[n{-}1]\mid \Des(uv_\setS)\subseteq\setS 
        \text{ and $\#\setS$ is odd}\}\,.
 \end{multline*}
\end{prop}

\begin{proof} 
 The formula above is equivalent to
 \begin{equation}\label{E:oneref}
   \lambda(u,v^{-1})\ =\ \lambda(v,u^{-1})\,,
 \end{equation}
 where $\lambda(\cdot,\cdot)$ is the structure constant for the antipode with
 respect to the fundamental basis, as proven in
 Theorem~\ref{T:ant-fundamental}.
 But~\eqref{E:oneref} expresses that $\Theta$ preserves antipodes (on
 the fundamental basis and its dual). 
\end{proof}

We turn now to quasi-symmetric functions. 
The dual $\QSym^*$ of $\QSym$ is the Hopf algebra  of {\em non-commutative
symmetric functions} of Gelfand, et.~al.~\cite{GKal}.  
It is the free associative algebra with generators
$\{M^*_{\emptyset,n}\mid n\geq 0\}$.
This statement is dual to  formula~\eqref{E:copqsym} for the coproduct of
$\QSym$. 

Define numbers
 \begin{align*}
   b(\setS,\setT)&\ :=\ \#\{u\in\frakS_n\mid \Des(u)\subseteq\setS,\
    \Des(u^{-1})\subseteq\setT\}\,,\\
   c(\setS,\setT)&\ :=\ \#\{u\in\frakS_n\mid \Des(u)\subseteq\setS,\
    \Des(u^{-1})\supseteq\setT\}\,,\\
   d(\setS,\setT)&\ :=\ \#\{u\in\frakS_n\mid
    \Des(u)=\setS,\ \Des(u^{-1})=\setT\}\,.
 \end{align*}
Let $\Phi$ denote the composite
 \[
   \QSym^*\ \xrightarrow{\ \calD^*\ }\ (\SSym)^*\ \xrightarrow{\ \Theta\ }\
    \SSym\ \xrightarrow{\ \calD\ }\ \QSym\,.
 \]

\begin{prop} \label{P:qsymdual} 
  The morphism $\Phi\colon \QSym^*\to\QSym$ sends
 \[
    F_\setS^*\ \mapsto\ \sum_{\setT\in\calQ_n} d(\setS,\setT)F_\setT
    \text{\ \ and \ \ }
     M_\setS^*\ \mapsto\ \sum_{\setT\in\calQ_n} b(\setS,\setT)M_\setT\,,
 \]
 for $\setS\in\calQ_n$.
\end{prop}

\begin{proof}  
 Since $\calD(\calF_u)=F_{\Des(u)}$, the dual  map satisfies
 $\calD^*(F_\setS^*)=\sum_{\Des(u)=\setS}\calF_u^*$. 
  Also, Theorem \ref{T:map-monomial} dualizes to
 $\calD^*(M_\setS^*)=\calM_{\zeta_\setS}^*$. 
 The descriptions of the composite above follow now from those for $\Theta$
 in~\eqref{E:selffund} and~\eqref{E:selfmon}, plus  that
 $\theta(\zeta_\setS,\zeta_\setT)=b(\setS,\setT)$, which in turn 
 follows from \eqref{E:galoisdes}.
\end{proof}

We now use the fact that $\Phi\colon\QSym^*\to\QSym$ is a morphism of
Hopf algebras.  
The image of $\Phi$ is precisely the 
subalgebra of $\QSym$ consisting of {\em symmetric} functions.
Since $\QSym^*$ is generated by $\{M^*_{\emptyset,n}\mid n\geq 0\}$, its 
image $\Phi(\QSym^*)$ is generated by $\Phi(M^*_{\emptyset,n})$, for $n\geq 
0$. Observe that $b(\emptyset,\setT)=1$ for every $\setT\in\calQ_n$ as $1_n$ is
the only permutation $u$ in $\frakS_n$ with $\Des(u)\subseteq\emptyset$ and 
$\emptyset=\Des(1_n^{-1})\subseteq\setT$.
Thus
 \[
   \Phi(M^*_{\emptyset,n})\ =\ 
     \sum_{\setT\in\calQ_n} M_\setT\ =\ F_{\emptyset, n}\,.
 \]
Formula~\eqref{E:QSym-def} shows that $F_{\emptyset,n}$ is the
complete homogeneous symmetric function of degree $n$. These
 generate the algebra of symmetric functions~\cite{Mac,St99}.
Thus, $\Phi$ is the
\emph{abelianization} map from non-commutative to commutative
symmetric functions. 
We will not use this, but rather the explicit expression
of $\Phi$ of Proposition~\ref{P:qsymdual}.

Let $a^\setR(\setS,\setT)$ denote the structure constants of the product
of $\QSym$ with respect to its monomial basis. These integers are
combinatorially described by \eqref{E:prodqsym} or \eqref{E:defalphaqsym}.
The following analog of Theorem~\ref{T:theta-form} provides a recursion for
computing the numbers $b(\setS,\setT)$ in terms of the structure constants 
$a^\setR(\setS,\setT)$. 
We view $a^\setR$ and $b$ as matrices with entries indexed by subsets of $[n{-}1]$.

\begin{prop} 
 For any $\setS\subseteq[p{-}1]$, $\setT\subseteq[q{-}1]$, 
 and $\setR\subseteq[p{+}q{-}1]$,
 \begin{equation}\label{E:b-prod}
  (ba^\setR b)(\setS,\setT)\ =\ b(\setS\cup\{p\}\cup(p+\setT),\setR).
 \end{equation}
\end{prop}
 
\begin{proof}  
 The dual of the coproduct of $\QSym$~\eqref{E:copqsym} is 
 \[
    M_\setS^*\cdot M_\setT^*\ =\ M_{\setS\cup\{p\}\cup(p+\setT)}^*\,.
 \]
 Thus, the right hand side of \eqref{E:b-prod} is the coefficient of $M_\setR$
 in $\Phi(M_\setS^*\cdot M_\setT^*)$. 
 On the other hand, since $b(\setS,\setT)=b(\setT,\setS)$, we have
 \[
   (ba^\setR b)(\setS,\setT)\ =\ 
    \sum_{\setS',\setT'} b(\setS,\setS')a^\setR(\setS',\setT')b(\setT',\setT)
    \ =\ 
    \sum_{\setS',\setT'}b(\setS,\setS')b(\setT,\setT')a^\setR(\setS',\setT')\,.
 \] 
 Thus the left hand side of~\eqref{E:b-prod} is the coefficient of $M_\setR$ in
 $\Phi(M_\setS^*)\cdot\Phi(M_\setT^*)$. 
 Since $\Phi$ is a morphism of algebras,~\eqref{E:b-prod} holds. 
\end{proof}

Expressing that $\Theta$ preserves the antipode in terms of the 
fundamental basis and its dual gives a result of Foata and
Sch\"utzenberger~\cite{FS78}, which Gessel obtained in his original work on
quasi-symmetric functions by other means~\cite[Corollary 6]{Ges} 
(Equation ({\it iv}) in the following corollary).
For $\setS\subseteq [n{-}1]$, define 
 \begin{align*}
   \setS^c&\ =\{i\in[n{-}1]\mid i\notin\setS\}\\
   \widetilde{\setS}&\ =\ \{i\in[n{-}1]\mid n-i\in\setS\}\,.
 \end{align*}

\begin{coro}\label{C:gessel}
 For $\setS,\setT\subseteq [n{-}1]$, the numbers $d(\setS,\setT)$ satisfy 
\begin{itemize} 
 \item[({\it i})] $d(\setS,\setT)=d(\setT,\setS)$\,,
 \item[({\it ii})] $d(\setS,\setT)=d(\widetilde{\setS},\widetilde{\setT})$\,, 
 \item[({\it iii})] $d(\setS,\setT)=d(\setS^c,\setT^c)$\,,  and
 \item[({\it iv})] $d(\setS,\setT)=d(\widetilde{\setS},\setT)$. 
\end{itemize}
\end{coro}

\begin{proof}
 The symmetry ({\it i}) follows by considering the bijection $u\mapsto u^{-1}$. 
 Similarly, ({\it ii}) follows by considering the bijection 
 $u\mapsto \omega_nu\omega_n^{-1}$, where  $\omega_n=(n,\ldots2,1)$,
 as it is easy to see that $\Des(\omega_nu\omega_n^{-1})=\widetilde{\Des(u)}$.

 The antipode of $\QSym$ is~\cite[corollaire 4.20]{Malv}
 \[
   S(F_\setT)\ =\ (-1)^nF_{\widetilde{\setT}^c}\,.
 \]
 Since $\Phi$ preserves antipodes, its explicit description in 
 Proposition~\ref{P:qsymdual} implies that 
 $d(\widetilde{\setS}^c,\setT)= d(\setS,\widetilde{\setT}^c)$. 
 Together with ({\it ii}) this yields ({\it iii}).

 Finally, to deduce ({\it iv}), consider the bijection $u\mapsto\omega_nu$. 
 Note that $\Des(\omega_nu)=\Des(u)^c$. 
 Therefore
 \[
   \Des((\omega_nu)^{-1})\ =\ 
   \Des(\omega_n\omega_nu^{-1}\omega_n^{-1})\ =\ 
   \Des(\omega_nu\omega_n^{-1})^c\ =\ 
   \widetilde{\Des(u^{-1})}^c\,.
 \]
 This shows that $d(\setS,\setT)=d(\setS^c,\widetilde{\setT}^c)$.
 Together with ({\it ii}) and ({\it iii}) this gives ({\it iv}).
\end{proof}

Expressing the preservation of the antipode under $\Phi$ in terms of monomial
quasi-symmetric functions and their duals gives further, similar results.

\begin{prop}\label{P:qsymdualantipode}
   The map  $S\Phi=\Phi S^*\colon \QSym^*\to\QSym$ sends
 \[
   M_\setS^*\ \mapsto\ (-1)^n\sum_\setR c(\setS,\widetilde{\setR}^c)M_\setR
    \ =\ (-1)^n\sum_\setR c(\setR,\widetilde{\setS}^c)M_\setR\,.
 \]
 Therefore,
 \[
    c(\setS,\widetilde{\setR}^c)\ =\ c(\setR,\widetilde{\setS}^c)\,.
 \]
\end{prop}

\begin{proof}
 We will show that 
 $S\Phi(M_\setS^*)=(-1)^n\sum_\setR c(\setS,\widetilde{\setR}^c)M_\setR$.
 One shows similarly that 
 $\Phi S^*(M_\setS^*)=(-1)^n\sum_\setR c(\setR,\widetilde{\setS}^c)M_\setR$.

 As mentioned in~\eqref{E:Q-antipode}, the antipode of $\QSym$ is
  \[
    S(M_\setT)\ =\ 
     (-1)^{\#\setT+1}\sum_{\setR\subseteq\setT}M_{\widetilde{\setR}}\,.
 \]
 Combining this with Proposition \ref{P:qsymdual} shows that $S\Phi$ sends
 \[
   M_\setS^*\ \mapsto\ 
    \sum_\setT b(\setS,\setT)(-1)^{\#\setT+1}
        \sum_{\setR\subseteq\setT}M_{\widetilde{\setR}}\,.
 \]
 Thus, we have to show that for each $\setS$ and $\setR$,
 \[
   \sum_{\setR\subseteq\setT}(-1)^{\#\setT+1} b(\setS,\setT)\ =\ 
    (-1)^n c(\setS,\setR^c)\,.
 \]
 Now,
 \begin{align*}
  \sum_{\setR\subseteq\setT}(-1)^{\#\setT+1} b(\setS,\setT) 
   &\ =\ 
    \sum_{\setR\subseteq\setT}\sum_{\setT'\subseteq\setT} (-1)^{\#\setT+1}
    \#\{u\mid \Des(u)\subseteq\setS,\ \Des(u^{-1})=\setT'\}\\
   &\ =\ \sum_{\setT'} \#\{u\mid \Des(u)\subseteq\setS,\ \Des(u^{-1})=\setT'\}
    \sum_{\setR\cup\setT'\subseteq\setT} (-1)^{\#\setT+1}\\
   &\ =\ \sum_{\setT':\setR\cup\setT'=[n{-}1]} (-1)^{n} 
     \#\{u\mid \Des(u)\subseteq\setS,\ \Des(u^{-1})=\setT'\}\\ 
   &\ =\  (-1)^{n} \#\{u\mid \Des(u)\subseteq\setS,\
    \Des(u^{-1})\cup\setR=[n{-}1]\}\\ 
   &\ =\ (-1)^{n} \#\{u\mid \Des(u)\subseteq\setS,\
        \Des(u^{-1})\supseteq\setR^c\} \\ 
   &\ =\ (-1)^{n}c(\setS,\setR^c)\,.    
 \end{align*}
\end{proof}

For completeness, we include the consequences on the numbers $b$ and $c$ that
follow. 
Note that these  also follow directly from 
 \[
   b(\setS,\setT)\ =\ 
     \sumsub{S'\subseteq\setS\\T'\subseteq\setT}d(\setS,\setT) 
    \ \text{ and }\ 
   c(\setS,\setT)\ =\ 
     \sumsub{S'\subseteq\setS\\T'\supseteq\setT}d(\setS,\setT)\,.
\]

\begin{coro}\label{C:gessel-monomial} 
 For any $\setS$, $\setT\subseteq[n{-}1]$,
\begin{itemize}
 \item[({\it i})]
            $b(\setS,\setT)=b(\setT,\setS)$,
 \item[({\it ii})]
            $b(\setS,\setT)=b(\widetilde{\setS},\widetilde{\setT})$, and 
              $c(\setS,\setT)=c(\widetilde{\setS},\widetilde{\setT})$,
 \item[({\it iii})]
            $c(\setS,\setT)=c(\setT^c,\setS^c)$,
 \item[({\it iv})]
            $b(\setS,\setT)=b(\widetilde{\setS},\setT)$, and 
              $c(\setS,\setT)=c(\widetilde{\setS},\setT)$.
\end{itemize}
\end{coro}

\def\cprime{$'$}
\providecommand{\bysame}{\leavevmode\hbox to3em{\hrulefill}\thinspace}
\providecommand{\MR}{\relax\ifhmode\unskip\space\fi MR }
\providecommand{\MRhref}[2]{%
  \href{http://www.ams.org/mathscinet-getitem?mr=#1}{#2}
}
\providecommand{\href}[2]{#2}


\begin{thebibliography}{10}

\bibitem{AF00}
Marcelo Aguiar and Walter Ferrer~Santos, \emph{Galois connections for incidence
  {H}opf algebras of partially ordered sets}, Adv. Math. \textbf{151} (2000),
  no.~1, 71--100. \MR{2001f:06007}

\bibitem{BilLiu}
Louis~J. Billera and Niandong Liu, \emph{Noncommutative enumeration in graded
  posets}, J. Algebraic Combin. \textbf{12} (2000), no.~1, 7--24.
  \MR{2001h:05009}

\bibitem{Bj}
Anders Bj{\"o}rner, \emph{Orderings of {C}oxeter groups}, Combinatorics and
  algebra (Boulder, Colo., 1983), Amer. Math. Soc., Providence, RI, 1984,
  pp.~175--195. \MR{86i:05024}

\bibitem{BLSWZ}
Anders Bj{\"o}rner, Michel Las~Vergnas, Bernd Sturmfels, Neil White, and
  G{\"u}nter~M. Ziegler, \emph{Oriented matroids}, second ed., Cambridge
  University Press, Cambridge, 1999. \MR{2000j:52016}

\bibitem{BCM}
Robert J. Blattner, Miriam Cohen and Susan Montgomery,
\emph{Crossed products and inner actions of Hopf algebras},
Trans. Amer. Math. Soc. \textbf{298} (1986), no.~2, 671--711.
\MR{87k:16012}

\bibitem{Co74}
Louis Comtet, \emph{Advanced combinatorics}, enlarged ed., D. Reidel Publishing
  Co., Dordrecht, 1974, The art of finite and infinite expansions. \MR{57
  \#124}

\bibitem{DHT}
G{\'e}rard Duchamp, Florent Hivert, and Jean-Yves Thibon, \emph{Some
  generalizations of quasi-symmetric functions and noncommutative symmetric
  functions}, Formal power series and algebraic combinatorics (Moscow, 2000),
  Springer, Berlin, 2000, pp.~170--178. \MR{2001k:05202}

\bibitem{DHT01}
G{\'e}rard Duchamp, Florent Hivert, and Jean-Yves Thibon,
\emph{Noncommutative symmetric functions VI: free quasi-symmetric functions
and related algebras}, Internat. J. Algebra Comput.  \textbf{12}  
(2002), no. 5, 671--717. \MR{1935570}

\bibitem{Ede} Paul Edeleman, \emph{Geometry and the
{M}\"obius function of the weak {B}ruhat   order of the symmetric group}, 1983.

\bibitem{Eh96}
Richard Ehrenborg, \emph{On posets and {H}opf algebras}, Adv. Math.
  \textbf{119} (1996), no.~1, 1--25. \MR{97e:16079}

\bibitem{FS78}
Dominique Foata and Marcel-Paul Sch{\"u}tzenberger, \emph{Major index and
  inversion number of permutations}, Math. Nachr. \textbf{83} (1978), 143--159.
  \MR{81d:05007}

\bibitem{GKal}
Israel~M. Gelfand, Daniel Krob, Alain Lascoux, Bernard Leclerc, Vladimir~S.
  Retakh, and Jean-Yves Thibon, \emph{Noncommutative symmetric functions}, Adv.
  Math. \textbf{112} (1995), no.~2, 218--348. \MR{96e:05175}

\bibitem{Ges}
Ira~M. Gessel, \emph{Multipartite ${P}$-partitions and inner products of skew
  {S}chur functions}, Combinatorics and algebra (Boulder, Colo., 1983)
  (Providence, RI), Amer. Math. Soc., 1984, pp.~289--317. \MR{86k:05007}

\bibitem{GR60a}
G.~Th. Guilbaud and P.~Rosenstiehl, \emph{Analyse alg\'ebrique d'un scrutin},
  M. Sci. humaines \textbf{4} (1960), 9--33.

\bibitem{Ho00}
Michael~E. Hoffman, \emph{Quasi-shuffle products}, J. Algebraic Combin.
  \textbf{11} (2000), no.~1, 49--68. \MR{2001f:05157}

\bibitem{Jol99}
Armin J\"ollenbeck, 
\emph{Nichtkommutative Charaktertheorie der symmetrischen Gruppen},
 [Noncommutative theory of characters of symmetric groups]
Bayreuth. Math. Schr. No. 56, (1999), 1--41.

\bibitem{Len72}
Andr\'e Lentin, \emph{\'Equations dans les mono\"\i des libres},
Math\'ematiques et Sciences de l'Homme, 
No. 16. Mouton; Gauthier-Villars, Paris, 1972. iv+160 pp. 
 
\bibitem{Lod}
Jean-Louis Loday, \emph{Homotopical syzygies}
Une d\'egustation topologique [Topological morsels]: homotopy theory in the
Swiss Alps (Arolla, 1999), Contemp. Math. \textbf{265} (2000), 99--127.
\MR{2001k:55022}

\bibitem{LR98}
Jean-Louis Loday and Mar{\'\i}a~O. Ronco, \emph{Hopf algebra of the planar
  binary trees}, Adv. Math. \textbf{139} (1998), no.~2, 293--309.
  \MR{99m:16063}

\bibitem{LR01}
\bysame, \emph{Order structure on the algebra of permutations and of planar
  binary trees}, 
J. Algebraic Combin. \textbf{15} (2002), no.~3, 253--270.
\MR{1900627}

\bibitem{Mac}
Ian G.~Macdonald, \emph{Symmetric functions and {H}all polynomials}, second ed.,
  The Clarendon Press Oxford University Press, New York, 1995, With
  contributions by A. Zelevinsky, Oxford Science Publications. \MR{96h:05207}

\bibitem{Malv}
Claudia Malvenuto, \emph{Produits et coproduits des fonctions
quasi-sym\'etriques et   de l'alg\`ebre des descents}, no.~16, Laboratoire de
combinatoire et   d'informatique math\'ematique {(LACIM)}, Univ.~du Qu\'ebec \`a
Montr\'eal,   Montr\'eal, 1994.

\bibitem{MR95}
Claudia Malvenuto and Christophe Reutenauer, \emph{Duality between
  quasi-symmetric functions and the {S}olomon descent algebra}, J. Algebra
  \textbf{177} (1995), no.~3, 967--982. \MR{97d:05277}

\bibitem{Mi66}
  R. James Milgram, \emph{Iterated loop spaces},
  Ann. of Math. (2) \textbf{84} (1966), 386--403.  \MR{34 \#6767}

\bibitem{MM65}
John~W. Milnor and John~C. Moore, \emph{On the structure of {H}opf algebras},
  Ann. of Math. (2) \textbf{81} (1965), 211--264. \MR{30 \#4259}

\bibitem{Mo93a}
Susan Montgomery, \emph{Hopf algebras and their actions on rings}, Published
  for the Conference Board of the Mathematical Sciences, Washington, DC, 1993.
  \MR{94i:16019}

\bibitem{RP01}
 Fr\'ed\'eric Patras and Christophe Reutenauer, \emph{Lie representations and an
  algebra containing {S}olomon's}, 
J. Algebraic Combin.  \textbf{16} (2002), no.~3, 301--314 (2003).

\bibitem{PR95}
St\'ephane Poirier and Christophe Reutenauer, \emph{Alg\`ebres de {H}opf de
  tableaux}, Ann. Sci. Math. Qu\'ebec \textbf{19} (1995), no.~1, 79--90.
  \MR{96g:05146}

\bibitem{Re93}
Christophe Reutenauer, \emph{Free {L}ie algebras}, The Clarendon Press Oxford
  University Press, New York, 1993, Oxford Science Publications. \MR{94j:17002}

\bibitem{Rot}
Gian-Carlo Rota, \emph{On the foundations of combinatorial theory. {I}.
  {T}heory of {M}\"obius functions}, Z. Wahrscheinlichkeitstheorie und Verw.
  Gebiete \textbf{2} (1964), 340--368 (1964), Reprinted in {\em Gian-Carlo Rota
  on Combinatorics: Introductory papers and commentaries} (Joseph P.S. Kung,
  Ed.), Birkh\"auser, Boston, 1995. \MR{30 \#4688}

\bibitem{SlSeek}
Neil J.~A.~Sloane, \emph{An on-line version of the encyclopedia of integer
  sequences}, Electron. J. Combin. \textbf{1} (1994), Feature 1, approx.\ 5
  pp.\ (electronic), {\tt
  http://akpublic.research.att.com/\~{}njas/sequences/ol.html}. \MR{95b:05001}

\bibitem{St86}
Richard~P. Stanley, \emph{Enumerative combinatorics. {V}ol. 1}, Cambridge
  University Press, Cambridge, 1997, With a foreword by Gian-Carlo Rota,
  Corrected reprint of the 1986 original. \MR{98a:05001}

\bibitem{St99}
\bysame, \emph{Enumerative combinatorics. {V}ol. 2}, Cambridge University
  Press, Cambridge, 1999, With a foreword by Gian-Carlo Rota and appendix 1 by
  Sergey Fomin. \MR{2000k:05026}

\bibitem{Stem97}
John~R. Stembridge, \emph{Enriched ${P}$-partitions}, Trans. Amer. Math. Soc.
  \textbf{349} (1997), no.~2, 763--788. \MR{97f:06006}

 \bibitem{Tak71}
Mitsuhiro Takeuchi, 
\emph{Free Hopf algebras generated by coalgebras},
J. Math. Soc. Japan \textbf{23} (1971), 561--582.

\bibitem{TU96}
Jean-Yves Thibon and B.-C.-V. Ung, \emph{Quantum quasi-symmetric functions 
  and Hecke algebras}, J. Phys. A \textbf{29} (1996), no.~22, 7337--7348. 
 \MR{97k:05204}



\end{thebibliography}
\end{document}